\documentclass[reqno]{amsart}
\usepackage{amsaddr}
\usepackage[english]{babel}

\usepackage[utf8]{inputenc}
\usepackage[T1]{fontenc}
\usepackage{lmodern}
\usepackage{microtype}
\usepackage{bm}
\usepackage{fullpage}

\usepackage{amsmath}
\usepackage{amsfonts}
\usepackage{amssymb}
\usepackage{amsthm}
\usepackage{mathtools}

\usepackage{graphicx}
\graphicspath{{figures/}}
\usepackage{booktabs}
\usepackage{algorithm}
\usepackage{algpseudocode}

\usepackage{xcolor}
\usepackage[colorlinks=true,
            linkcolor=blue!55!black,
            citecolor=blue!55!black,
            urlcolor=blue!55!black,
            linktoc=all]{hyperref}
\usepackage{orcidlink}

\newtheorem{theorem}{Theorem}[section]
\newtheorem{lemma}[theorem]{Lemma}
\newtheorem{assumption}[theorem]{Assumption}
\newtheorem{proposition}[theorem]{Proposition}

\theoremstyle{remark}
\newtheorem{remark}[theorem]{Remark}
\newcommand{\bydefeq}{\coloneqq} 
\let\origmarginpar\marginpar
\renewcommand{\marginpar}[1]{\origmarginpar{\scriptsize #1}} 
\newcommand{\rd}{\mathop{}\!\mathrm{d}}
\newcommand{\PP}{\mathbb{P}}
\newcommand{\R}{\mathbb{R}}
\newcommand{\N}{\mathbb{N}}
\newcommand{\Rmn}{\mathbb{R}^{m\times n}}
\newcommand{\Rm}{\mathbb{R}^{m}}
\newcommand{\Rs}{\mathbb{R}^{s}}
\newcommand{\Rn}{\mathbb{R}^{n}}
\newcommand{\tr}{^\top}
\newcommand{\zb}{\bm{z}}
\newcommand{\fb}{\bm{f}}
\newcommand{\yb}{\bm{y}}
\newcommand{\kappaeff}{\kappa_{\text{eff}}}

\newcommand{\vb}{\bm{v}}
\newcommand{\cMr}{\mathcal{M}_r}
\newcommand{\cP}{\mathcal{P}}
\newcommand{\cL}{\mathcal{L}}
\newcommand{\cJ}{\mathcal{J}}
\newcommand{\trace}{\mathrm{tr}}
\newcommand{\Jct}{\widetilde{\mathcal{J}}}
\newcommand{\Sct}{\widetilde{\mathcal{S}}}
\newcommand{\op}{\text{op}}

\usepackage{accents}
\newcommand*{\dt}[1]{%
 \accentset{\mbox{\large\bfseries .}}{#1}}
\newcommand*{\ddt}[1]{%
 \accentset{\mbox{\large\bfseries .\hspace{-0.25ex}.}}{#1}}

\allowdisplaybreaks
\author[B.~Carrel and T.~Vanzan]{Benjamin Carrel$^{{\lowercase\mathrm{a}}}$\orcidlink{0000-0003-3549-2845},
                \and 
		Tommaso Vanzan$^{{\lowercase\mathrm{b}}}$\orcidlink{0000-0001-7554-4692}}

\address[Benjamin Carrel]{
	$\,^{\lowercase\mathrm{a}}$ Independent researcher
	} 
\address[Tommaso Vanzan]{
	$\,^{\lowercase\mathrm{b}}$Dipartimento di Scienze Matematiche, 
	Politecnico di Torino \\
  	Corso Duca degli Abruzzi 24, 10129 Torino,  
	Italy}

\title[DLRA for the optimal control of matrix ODEs]{Dynamical low-rank approximations for the optimal control of matrix ODEs}
\date{\today}

\keywords{Dynamical low-rank approximation, Optimal control, Parabolic random PDEs, Covariant ODEs, Matrix differential equations}

\subjclass[2020]{
			 65K10,  	
			 65L20, 
			  35R60, 
			  65N22.
			    } 

\begin{document}

\begin{abstract}
We consider the open-loop optimal control of a dynamical system whose state is described by a time-dependent matrix $X(t)\in\Rmn$. Standard gradient-based optimization methods may incur substantial computational and memory costs, particularly for large-scale systems and fine temporal discretizations. To alleviate this computational burden, we develop a mathematical framework to use Dynamical Low-Rank Approximations (DLRA) in optimal control problems. We study two distinct paradigms: the first one, called \textit{Optimize-then-Project} (OtP), first derives the optimality system and then approximates both the state and adjoint equations using DLRA. The second one, called \textit{Project-then-Optimize} (PtO), first replaces the state equation with its DLR approximation, and then derives the optimality conditions for the resulting optimization problem. The two approaches differ fundamentally concerning the role of the adjoint variable. In OtP, the adjoint variable evolves on the manifold of fixed-rank matrices, whereas in PtO, the adjoint variable belongs to the tangent space to the fixed-rank manifold at the evolving state variable. This latter formulation leads to a covariant Ordinary Differential Equation (ODE) describing the dynamics of the adjoint variable as the underlying tangent space evolves in time.
 
For both approaches we establish rigorous error estimates for the approximated gradients, identifying how their accuracy depends on the low-rank approximations of initial and terminal conditions, the components of the state and adjoint vector fields normal to the relevant tangent spaces, the stability properties of the dynamics, and the manifold curvature.
The PtO approach further requires tailored time discretization strategies for the resulting adjoint covariant ODE,
and we propose one with an error bound that depends on the variation of the tangent projector along the state trajectory.
To illustrate the proposed framework, we apply these methodologies to the optimal control of a parabolic partial differential equation with random coefficients and initial conditions.
\end{abstract}

\maketitle	    
\section{Introduction}\label{sec:intro}
In this manuscript, we address the efficient numerical solution of open-loop optimal control problems constrained by matrix Ordinary Differential Equations (ODEs). The prototype example reads
\begin{equation}\label{eq:intro_modelproblem}
\begin{aligned}
  & \min_{X,\; \zb} \mathcal{J}(X, \zb)
  = \tfrac{1}{2} \int_0^T \!\! \| X(t) - X_d(t) \|^2 \, \rd t
      \;+\; \tfrac{\nu}{2} \int_0^T \!\! \| \zb(t) \|^2 \, \rd t,\\
      &\text{subject to }\; \begin{cases}
  \dt X(t) &=\; F\bigl(t, X(t),\, \zb(t) \bigr),\quad t\in (0,T],\\
   X(0) &= X_0.
\end{cases}
\end{aligned}
\end{equation}
with $X(t)\in \Rmn$ and $\zb(t)\in \Rs$, $t\in [0,T]$. Many dynamical systems, often formulated in vectorized form, can be expressed as matrix (or more generally tensor) ODEs. Examples include time-dependent second-order partial differential equations (PDEs) semidiscretized spatially on Cartesian meshes \cite{DAUTILIA20202067}, density matrix equations in quantum mechanics \cite{blum2012density}, kinetic equations \cite{einkemmer2018low} and random PDEs \cite{sapsis2009dynamically}.
To compute a minimizer of \eqref{eq:intro_modelproblem}, a standard gradient-based iteration requires a forward solve for the state matrix ODE and a backward solve for the adjoint matrix ODE, followed by an update of the control. When the time interval is discretized with $N_t$ subintervals, and the state and adjoint variables are stored as dense $m\times n$ matrices at every time step, the memory footprint and computational cost typically scale as $\mathcal{O}(mnN_t)$, which can severely limit numerical computations in practical scenarios.
A way to reduce this cost is to exploit the potential low-rank approximability of the state and adjoint matrices.

Dynamical low-rank approximation (DLRA) \cite{kochlubich2007,nonnenmacher2008dynamical} has emerged as a powerful framework to compute low-rank approximations of time-dependent matrices. Specifically, DLRA approximates $X(t)$ with a low-rank matrix $Y(t)=U(t)S(t)V(t)^\top\in \Rmn$ that belongs to the manifold $\cMr$ of matrices of fixed rank $r\in\mathbb{N}$. The time evolution for $Y(t)$ is derived by imposing that, at each time instant, the time derivative $\dt Y(t)$ is equal to the orthogonal projection of the full vector field $F$ onto the tangent space at $Y(t)$. This construction yields a reduced dynamical system that evolves only the low-rank factors rather than the full dense matrix. Consequently, DLRA reduces the memory footprint at each time step from $\mathcal{O}(mn)$ to $\mathcal{O}((m+n)r)$, and leads also to comparable savings in computational times if the vector field can be evaluated directly from the low-rank factors.
After several contributions aimed at developing time integrators for the resulting differential equations that were robust against small singular values of $Y(t)$ \cite{lubichoseledets2014,kieri2016robust,kieri2019projection,ceruti2022unconventional,ceruti2022rank}, DLRA has found numerous applications in a variety of fields, ranging from kinetic and plasma equations \cite{einkemmer2018low,einkemmer2025review}, random PDEs \cite{musharbash2020symplectic,sapsis2009dynamically}, stochastic differential equations \cite{kazashi2025dynamical,bao2026exponential} and data assimilation \cite{schmidt2023rank,TRIGOTRINDADE2027115433,marzouk2026dynamical}.

Specifically in the context of optimization and control,
\cite{baumann2025adaptive} considers a parameter identification problem for kinetic equations and proposes an ``optimize then discretize then low-rank'' strategy. The work \cite{einkemmer2026memory} develops a memory-efficient adjoint optimization that avoids the storage of the low-rank state factors by recomputing them on-the-fly during the backward adjoint integration using a tailored time-reversible integrator for the state equation. Concerning closed-loop optimal control, \cite{Saluzzi2025} introduces a DLRA strategy for state-dependent Riccati equations for parametric dynamical systems, while \cite{eigel2023hjb} discusses a tensor-train approximation to solve the Hamilton--Jacobi--Bellman equation. DLRA has also been investigated for neural network training in \cite{schotthofer2022low,kusch2026augmented}.
Further low-rank strategies in optimization, but not DLRA-based, are discussed in \cite{BENNER201626,garreis2017constrained,benner2020,goutte2026analytic}.

Despite the substantial computational speed-ups reported in the above-mentioned works, a systematic mathematical framework for incorporating DLRA into optimal control is still lacking. The goal of this manuscript is to consider the prototype problem \eqref{eq:intro_modelproblem} and establish mathematical foundations for a DLRA-based optimization.

\subsection{Contributions and outline.}
Our contributions are fourfold.

First, we formulate and analyze two DLRA-based approaches for computing approximate gradients. The first one, called \textit{Optimized-then-Project} (OtP), derives the optimality system of \eqref{eq:intro_modelproblem} and then projects the vector fields of both the forward and adjoint equations onto the tangent spaces at the current state and adjoint variables. In this way, one derives two independent DLRAs, one for the state variable and one for the adjoint variable.
In contrast, the second approach, called \textit{Project-then-Optimize} (PtO), first replaces the state equation with its DLRA and then derives the optimality conditions of the approximate problem. In this case, the adjoint variable is not a point on $\cMr$; instead it is a vector belonging to the tangent space to $\cMr$ at the current state. Since the tangent space evolves as the state changes, the adjoint equation is expressed as a covariant ODE, which prescribes only the projected time derivative of the adjoint variable. 
Consequently, the two approaches are not equivalent, and the second one is particularly rich from a mathematical point of view.

Second, we establish rigorous error estimates for the gradient approximations yielded by both approaches.
The estimates specify how the accuracy of the approximated gradient depends on the low-rank approximation of the initial state, on the normal components of the state and adjoint vector fields to the relevant tangent spaces, and on the stability properties of the dynamics. For the OtP approach, the error bound also accounts for an additional term due to the noncoherent initialization to the zero matrix of the adjoint variable at the terminal time $T$ (as the zero matrix does not belong to $\cMr$). For the PtO paradigm, the error bound has instead two additional contributions involving the effective curvature, measuring the variation of the orthogonal projector as the state variable evolves, and a normal component of the state vector field.
Both estimates are consistent in the full-rank limit, meaning that when $r=\min\{m,n\}$, both OtP and PtO approaches recover the full-rank gradient.

Third, we discuss time integration strategies for both the forward and backward equations, with a particular focus on the backward covariant ODE of the PtO approach. We derive differential equations for the factors of the unique tangent space parametrization and argue that they are in general not stiff even when the smallest singular value of $Y(t)$ tends to zero, but they do depend on the \textit{effective curvature}, that is, on the actual speed of rotation of the factors $U(t)$ and $V(t)$. The use of these differential equations is however limited, since they require the knowledge of both $Y(t)$ and of its time derivative at every instant of time, which is generally unavailable when $Y(t)$ is itself a numerical solution computed only at a finite number of time instants. Therefore, we propose and analyze a projected Euler method that advances the adjoint variable by projecting the adjoint vector field onto the tangent space at the next state. This method represents a local second-order approximation to the differential equations for the tangent factors, and we prove its first-order convergence with an error constant that depends on the effective curvature, but not explicitly on the smallest singular value of $Y(t)$.

Fourth, we apply the proposed methodologies to the optimal control of a random parabolic PDE, with random coefficients and initial data. While the optimal control of stationary random PDEs has been the subject of intense research in the last decade \cite{kouri2016risk,kouri2018optimization,heinkenschloss2025optimization,nobile2025multilevel}, very few works have dealt with time-dependent equations, see, e.g., \cite{guth2024parabolic}, due to the extremely large computational burden.
For this application, we present numerical experiments aimed at verifying the theoretical assumptions and the error bounds, at testing how the approximations affect the minimization of the desired functional, and at assessing the computational speed-ups.

\subsection{Layout}
The rest of the manuscript is organized as follows. Subsection \ref{subsec:Notation} introduces the notation used in the manuscript. Section \ref{sec:setup} presents the model problem and discusses its well-posedness and optimality conditions. Section \ref{sec:manifold_DLRA} recalls the concepts on manifolds and on dynamical low-rank approximations that will be relevant for this work. Section \ref{sec:DLRAOPT} is the core of the manuscript and discusses in Subsection \ref{subsec:OtP} the optimize-then-project, together with the derivation of error estimates for the gradient approximation. Subsection \ref{subsec:PtO} is concerned with the project-then-optimize paradigm and contains the analysis of the related optimal control problem, the derivation of the optimality condition and of the error estimates. Section \ref{sec:disc} discusses time integration strategies for the forward and backward equations. Subsection \ref{subsubsec:covariantODE} is concerned with numerical methods to integrate covariant ODEs. Section \ref{sec:ocp} presents an application to the optimal control of random parabolic equations, together with supporting numerical experiments. Finally, Appendix \ref{appendix:adjoint} details the derivation of the adjoint of the derivative of the projection operator.

\subsection{Notation}\label{subsec:Notation}
Let $M\in \mathbb{R}^{m\times m}$, $\Xi \in \mathbb{R}^{n\times n}$ and $\Gamma\in \mathbb{R}^{s\times s}$ be symmetric positive definite. We introduce the weighted inner products,
\begin{align}
\langle X, Y\rangle_{M, \Xi}&\bydefeq \trace\bigl( X^\top M Y\, \Xi \bigr),\; \forall X,\; Y\in \mathbb{R}^{m\times n},\quad  \langle \mathbf{u},\vb \rangle_{M} \bydefeq \mathbf{u}^\top M \vb,\; \forall \mathbf{u},\; \vb\in \mathbb{R}^{m},\\
\langle \mathbf{u},\vb \rangle_{\Xi} &\bydefeq \mathbf{u}^\top \Xi \vb,\;\forall \mathbf{u},\; \vb\in \mathbb{R}^{n},\hspace*{2.2cm} \langle \mathbf{u},\vb \rangle_{\Gamma} \bydefeq \mathbf{u}^\top \Gamma \vb,\;\forall \mathbf{u},\; \vb\in \mathbb{R}^{s},
\end{align}
respectively, and denote the associated norms by
$\|\cdot\|_{M,\Xi}$, $\|\cdot\|_M$, $\|\cdot\|_\Xi$, and
$\|\cdot\|_\Gamma$.
Throughout the manuscript, we always consider the spaces $\mathbb{R}^{m\times n}$, $\mathbb{R}^m$, $\mathbb{R}^n$ and $\Rs$ equipped with these specific inner products and induced norms.
Let $I$ be an open interval of $\R$ and $p\in [1,\infty)$. We denote by $L^p(I; \Rs)$ the Bochner space of measurable functions such that $\|\phi\|_{L^p(I;\Rs)}=\left(\int_I \|\phi(t)\|_{\Gamma}^p\rd t\right)^{\frac{1}{p}}<\infty$, and by $H^k(I;\Rs)$, $k\in\N$, the $k$-th order Sobolev space of $L^2$-functions.
The space of all linear functionals from $\Rm$ to $\Rm$ is denoted by $\cL(\Rm)$, while the space of linear maps from $\Rs$ to $\Rmn$ is denoted by $\cL(\Rs;\Rmn)$ and equipped with the norm $\|B\|_{\op}:=\sup\{\|B\vb\|_{M,\Xi}:\;\vb\in \Rs:\;\|\vb\|_\Gamma \leq 1\}$.

\section{Problem setting}
\label{sec:setup}
Let $T>0$ and $I:=(0,T)$ be a finite time interval. In this work we are concerned with the efficient solution of the optimal control problem
  \begin{align}
  & \min_{X\in H^1(I;\mathbb{R}^{m\times n}),\; \zb\in L^2(I;\Rs)} \mathcal{J}(X, \zb)
  \bydefeq \tfrac{1}{2} \int_0^T \!\! \| X(t) - X_d(t) \|_{M, \Xi}^2 \, \rd t
      \;+\; \tfrac{\nu}{2} \int_0^T \!\! \| \zb(t) \|_{\Gamma}^2 \, \rd t,\label{eq:ocp_full}\\
      &\text{subject to }\; \begin{cases}
  \dt X(t) &=\; F\bigl(t, X(t),\, \zb(t) \bigr),\quad \text{for a.e. }t\in (0,T],\\
   X(0) &= X_0.
  \end{cases}\label{eq:state_equation}
  \end{align}
The objective functional consists of a tracking term, which measures the discrepancy between the state and a desired trajectory, and of a Tikhonov regularization term weighted by the parameter $\nu>0$, which penalizes large control inputs.
In addition, $F:\R\times \Rmn\times \Rs\rightarrow \Rmn$ is a possibly nonlinear function describing the system dynamics, $X_d\in L^2(I;\Rmn)$ denotes the desired target trajectory, $X_0\in \Rmn$ is the prescribed initial condition, and $\zb\in L^2(I;\Rs)$ is the open-loop control.
We make the following assumption on the vector field $F$.

\begin{assumption}[Assumption on the vector field $F$]\label{ass:vector_field}
The vector field $F:\R\times \Rmn\times \Rs\rightarrow \Rmn$ takes the form 
\[F(t,X,\zb)=f(t,X)+B(t,X)\zb.\]
We assume that $f:\R\times \Rmn\rightarrow \Rmn$ and $B:\R\times \Rmn\rightarrow \cL(\Rs,\Rmn)$ are measurable in $t$ and globally Lipschitz continuous with respect to $X$, meaning that there exist two positive constants $L_f$ and $L_B$ such that for almost all $t\in I$ and all $X_1,X_2\in \Rmn$,
\begin{equation}\label{eq:Lipschitz}
\begin{aligned}
\|f(t,X_1)-f(t,X_2)\|_{M,\Xi} &\leq L_f\|X_1-X_2\|_{M,\Xi},\\
\|B(t,X_1)-B(t,X_2)\|_{\op} &\leq L_B\|X_1-X_2\|_{M,\Xi}.
\end{aligned}
\end{equation}
In addition, they satisfy the growth conditions $\|f(t,0)\|_{M,\Xi}\leq \alpha(t)$ and $\|B(t,0)\|_{\op}\leq \beta(t)$ for some $\alpha \in L^2(I;\R)$ and $\beta \in L^\infty(I;\R)$.
\end{assumption}
The function $f$ determines the uncontrolled system dynamics, while $B$ is the control operator. 
The hypotheses on $B$ are sufficiently general to include, e.g., bilinear controls of the form
\[B(t,X)\zb(t)=\sum_{i=1}^s z_i(t)B_i X,\quad B_i\in\mathbb{R}^{m\times m},\]
or controls that simultaneously act over all columns of $X$,
\[B(t,X)\zb(t)= (E_{s\rightarrow m}\zb(t))\mathbf{1}_n^\top,\]
as typically appearing in the robust optimal control of random PDEs (see Section \ref{sec:ocp}), where $E_{s\rightarrow m}:\Rs\rightarrow \Rm$ can be the extension operator of a boundary control and $\mathbf{1}_n\in \Rn$ is the vector containing all ones.

Assumption \ref{ass:vector_field} implies that for any $X\in \Rmn$,
\[\|F(t,X,\zb(t))\|_{M,\Xi}\leq\underbrace{\left(L_f+L_B\|\zb(t)\|_{\Gamma}\right)}_{=:l(t)}\|X\|_{M,\Xi}+\underbrace{\alpha(t)+\beta(t)\|\zb(t)\|_{\Gamma}}_{=:h(t)}=l(t)\|X\|_{M,\Xi}+h(t).\]
Since both $l$ and $h$ are in $L^1(I;\R)$, the Carath\'eodory existence theorem (see, e.g., \cite[Theorem 54 and Proposition C.3.8]{sontag2013mathematical} or \cite{walter2013ordinary}) guarantees that there exists a unique absolutely continuous function $X$ satisfying \eqref{eq:state_equation}. Furthermore, as $l$ and $h$ are even in $L^2(I;\mathbb{R})$, $\dt X \in L^2(I;\Rmn)$.
Consequently, we can define a solution operator $\mathcal{S}: L^2(I;\Rs)\rightarrow H^1(I;\Rmn)$ and introduce the reduced functional $\widehat{\mathcal{J}}(\zb)=\mathcal{J}(\mathcal{S}(\zb),\zb)$. The next lemma discusses the existence of a minimizer of the optimal control problem.

\begin{lemma}\label{lem:existence_minimizer_ocp}
Let Assumption \ref{ass:vector_field} hold.
Then, 
\begin{equation}\label{eq:ocp}
\min_{\zb \in L^2(I;\Rs)} \widehat{\mathcal{J}}(\zb)
\end{equation}
admits a minimizer.
\end{lemma}
\begin{proof}
Consider a minimizing sequence $\{\zb_k\}_{k\in \mathbb{N}}$.
Since $\nu>0$ and the tracking term is nonnegative, the map $\zb\mapsto \mathcal{J}(\mathcal{S}(\zb),\zb)$ is coercive, implying that there exists a $C_{\nu}$ such that $\|\zb_k\|_{L^2(I;\Rs)}\leq C_\nu$ for every $k\in\mathbb{N}$. Due to the reflexivity of $L^2(I;\Rs)$, there exists a subsequence, still denoted by $\{\zb_{k}\}_{k\in \mathbb{N}}$, converging weakly to $\zb^\star$ in $L^2(I;\Rs)$. 
Next, we show that $X_k\bydefeq\mathcal{S}(\zb_k)$ converges, up to subsequences, strongly to $\mathcal{S}(\zb^\star)$ in $C^0(\overline{I};\mathbb{R}^{m\times n})$.
To this end, observe that for all $\zb\in L^2(I;\Rs)$ and a.e. $t\in I$,
\begin{equation}\label{eq:bound_norm_derivative_X}
\begin{aligned}
\frac{d}{dt}\frac{1}{2}\| X(t)\|^2_{M,\Xi}=\langle X(t),\dt X(t)\rangle_{M,\Xi}&\leq \|X(t)\|_{M,\Xi}\|F(t,X(t),\zb(t))\|_{M,\Xi}\\
&\leq \left(L_f+L_B\|\zb(t)\|_\Gamma\right)\|X(t)\|^2_{M,\Xi}+\left(\alpha(t)+\beta(t)\|\zb(t)\|_\Gamma\right)\|X(t)\|_{M,\Xi}.\\
\end{aligned}
\end{equation}
Introducing the shorthand notation $E(t)\bydefeq\frac{1}{2}\|X(t)\|^2_{M,\Xi}$ and using Young's inequality, the previous relation can be written as
\[E'(t)\leq \underbrace{(1+2L_f+2L_B\|\zb(t)\|_\Gamma)}_{=:\mu(t)} E(t) +\underbrace{\frac{1}{2} \left(\alpha(t)+\beta(t)\|\zb(t)\|_\Gamma\right)^2}_{=:g(t)}\quad \text{for a.e. }t\in I,\]
which implies
\[E(t)\leq E(0)+\int_0^t g(s)\rd s + \int_0^t \mu(s) E(s)\rd s \quad\text{for every } t\in \overline{I}. \]
Gr\"onwall's lemma then leads to the upper bound
\begin{equation*}
E(t)\leq e^{\int_0^t \mu(s)\rd s} \left(E(0)+\int_0^t g(s)\rd s\right) \quad\text{for every } t\in \overline{I},
\end{equation*}
which allows to conclude that
\[\sup_{t\in[0,T]} \|X(t)\|^2_{M,\Xi}\leq 2 e^{T(1+2L_f)+2L_B\sqrt{T}\|\zb\|_{L^2(I;\Rs)}}\left(E(0)+\|\alpha\|^2_{L^2(I;\R)}+\|\beta\|^2_{L^\infty(I;\R)}\|\zb\|^2_{L^2(I;\Rs)}\right).\]
As the sequence $\{\zb_k\}_{k\in \N}$ is bounded, $\{X_k\}_{k\in \mathbb{N}}$ is uniformly bounded in the norm of $C^0(\overline{I};\mathbb{R}^{m\times n})$ by a constant $C_X=C_X(C_\nu)$.
Furthermore, the sequence $\{X_k\}_{k\in\mathbb{N}}$ is equicontinuous: for any $0\leq s\leq t\leq T$,
\begin{equation*}
\begin{aligned}
\|X_k(t)-X_k(s)\|_{M,\Xi}\leq \int_s^t \|F(\tau, X_k(\tau),\zb_k(\tau))\|_{M,\Xi}\rd \tau &\leq \int_s^t (l_k(\tau)\|X_k(\tau)\|_{M,\Xi}+h_k(\tau)) \rd \tau\\
&\leq \sqrt{t-s}\left(C_X C_l+C_m\right),
\end{aligned}
\end{equation*}
$C_l$ and $C_h$ being the uniform bounds in $k$ of $l_k(t):=L_f+L_B\|\zb_k(t)\|_{\Gamma}$ and $h_k(t):=\alpha(t)+\beta(t)\|\zb_k(t)\|_{\Gamma}$ in the norm of $L^2(I;\mathbb{R})$.
Due to the Arzel\`a--Ascoli theorem, we can extract a subsequence $\{X_k\}_{k\in \mathbb{N}}$ that converges strongly to some $X^\star$ in $C^0(\overline{I};\mathbb{R}^{m\times n})$.
Next, we identify $X^\star$ with $\mathcal{S}(\zb^\star)$. To do so, we consider the integral representation
\[X_k(t)=X_{0} +\int_0^t F(s,X_k(s),\zb_k(s)) \rd s=X_{0} +\int_0^t (f(s,X_k(s)) + B(s,X_k(s))\zb_k(s)) \rd s,\]
and discuss the limits of each term as $k\rightarrow \infty$.
On the one hand, due to the uniform convergence of $\{X_k\}_{k\in\mathbb{N}}$ to $X^\star$ and the Lipschitz continuity of $f$ with respect to $X$, we have that
$\int_0^t f(s,X_k(s))\rd s\rightarrow \int_0^t f(s,X^\star(s))\rd s$ in $\Rmn$ as $k\rightarrow \infty$.
On the other hand,
\begin{equation*}
\begin{split}
\left\|\int_0^t (B(s,X_k(s))\zb_k(s)-B(s,X^\star(s))\zb^\star(s)) \rd s\right\|_{M,\Xi} &\leq  \left\|\int_0^t \left(B(s,X_k(s))-B(s,X^\star(s))\right)\zb_k(s)\rd s\right\|_{M,\Xi}\\
& +\left\|\int_0^t B(s,X^\star(s))(\zb_k(s)-\zb^\star(s))\rd s\right\|_{M,\Xi}.
\end{split}
\end{equation*}
The first term tends to zero as $k\rightarrow\infty$ due to the Lipschitz continuity of $B$ with respect to $X$, the uniform convergence of $\{X_k\}_{k\in\mathbb{N}}$ and the uniform boundedness in $k$ of $\{\zb_k\}_{k\in\mathbb{N}}$. The second term tends to zero due to the weak convergence of $\{\zb_k\}_{k\in\mathbb{N}}$ to $\zb^\star$ in $L^2(I;\Rs)$.   
Hence, we conclude that $X^\star=\mathcal{S}(\zb^\star)$, and that $\mathcal{S}(\zb_k)$ converges to $\mathcal{S}(\zb^\star)$ in $C^0(\overline{I};\Rmn)$. 
Combining the strong continuity of the tracking term and the weak lower semicontinuity of the regularization term, we obtain 
\[\inf_{\zb\in L^2(I;\Rs)} \cJ(\mathcal{S}(\zb),\zb)=\liminf_{k\rightarrow \infty} \cJ(\mathcal{S}(\zb_k),\zb_k)\geq \cJ(\mathcal{S}(\zb^\star),\zb^\star),\]
and hence $\zb^\star$ is a minimizer of \eqref{eq:ocp}.
\end{proof}

Under the additional assumption that $f$ and $B$ are Fr\'echet differentiable with respect to $X$, every minimizer of \eqref{eq:ocp} satisfies the optimality system
\begin{align}
 \dt X(t) &= F(t,X(t),\zb(t))\quad & \text{for a.e. } t\in I,\label{eq:ocp_optimality_1}\\
 - \dt P(t) &= F_X^\star(t,X(t),\zb(t))P(t) + X_d(t)-X(t)\quad & \text{for a.e. }t\in I,\label{eq:ocp_optimality_2}\\
 \nu \zb(t) &= B^\star(t,X(t)) P(t)\quad & \text{for a.e. } t\in I,\label{eq:ocp_optimality_3}
\end{align}
equipped with the initial and terminal conditions $X(0)=X_0$ and $P(T)=0$. Here, $F_X^\star$ denotes the adjoint of the Fr\'echet derivative of $F$ with respect to $X$ at a point $(t,X,\zb)$, that is the operator satisfying
\[\langle Z, F_X^\star(t,X,\zb) W\rangle_{M,\Xi}=\langle F_X(t,X,\zb) Z,  W\rangle_{M,\Xi}\quad \forall Z,\;W\in \Rmn,\]
and a similar definition holds for $B^\star$. The well-posedness of the \textit{linear} adjoint equation can be shown by adapting the arguments used for \eqref{eq:state_equation}. In particular, the adjoint solution satisfies $P\in H^1(I;\mathbb R^{m\times n})$.

The minimization of $\widehat{\mathcal{J}}$ via a gradient descent algorithm would then entail the solution of the forward matrix equation \eqref{eq:ocp_optimality_1}, of the backward adjoint equation \eqref{eq:ocp_optimality_2}, and the update of the current iterate using \eqref{eq:ocp_optimality_3}. This procedure can be very computationally intensive whenever $m$ and $n$ are large.
The goal of this manuscript is to analyze algorithmic frameworks based on dynamical low-rank approximations to accelerate the computation of gradients. These are presented in Section \ref{sec:DLRAOPT}. To prepare the ground,
the next section recalls basic notions on the manifold of fixed-rank matrices and of the DLRA framework, together with a few technical derivations that will be needed in the following.

\section{The manifold of fixed-rank matrices and dynamical low-rank approximations}
\label{sec:manifold_DLRA}
Let $r\in \mathbb{N}$ and consider the set of real $m\times n$ matrices of rank $r$,
\begin{equation*}
  \label{eq:Mr}
  \cMr
  \;=\; \{\, X \in \mathbb{R}^{m \times n} \;\mid\; \mathrm{rank}(X) = r \,\}.
\end{equation*}
It is well-known that $\cMr$ is a smooth embedded submanifold of $\mathbb{R}^{m\times  n}$; see, e.g., \cite[Sec.~7.5]{boumal_2023}. In this work, we represent elements on $\cMr$ via the (nonunique) weighted parametrization\footnote{We employ weighted inner products as in many applications it can be desirable to use discrete inner products that are consistent discretizations of infinite-dimensional counterparts; for an example see Section \ref{sec:num} or \cite{sutti2026low}.}
\begin{equation}\label{eq:newparametrizationmanifold}
    \begin{split}
       \cMr = \lbrace US V\tr \colon & U\in \R^{m\times r},\; V\in \R^{n\times r},\;\text{s.t. }U\tr\! M U = V\tr\! \Xi V = I_{r}, \\
              & S \in\mathbb{R}^{r\times r},\;\text{$S$ invertible} \rbrace,
    \end{split}
\end{equation}
and elements in the tangent space to $\cMr$ at $Y = U S V^\top$, denoted by $T_Y \cMr$, via the unique parametrization
\begin{equation*}
\begin{aligned}
  T_Y \cMr
  \;=\; \bigl\{\,  U\, L\, V^\top
                \;+\; U_p\, V^\top
                \;+\; U\, V_p^\top
              \,&\bigm|\, U_p^\top M\, U = 0,\ V_p^\top \Xi\, V = 0,\\
                 & U_p\in\mathbb{R}^{m\times r},\; V_p\in\mathbb{R}^{n\times r},\;  L \in \mathbb{R}^{r \times r} \,\bigr\}.
\end{aligned}
\end{equation*}
The orthogonal complement to $T_{Y}\cMr$ in $\mathbb{R}^{m\times n}$ with respect to the $\langle\cdot,\cdot\rangle_{M,\Xi}$ inner product is called normal space and denoted by $N_{Y}\cMr$.
The orthogonal projector onto the tangent
space of $\cMr$ at $Y$ with respect to $\langle\cdot,\cdot\rangle_{M, \Xi}$ is denoted by $\cP_Y$ and its orthogonal complement projector by $\cP_Y^\perp$.
A direct calculation, see \cite[Lemma 4.1]{kochlubich2007}, shows that, for $Y=USV\tr \in \cMr$,
$\cP_Y$ admits the closed form
\begin{equation}\label{eq:proj-formula2}
\cP_Y Z = P_U Z + ZP_V\tr - P_U Z P_V\tr= Z - P_U^\perp Z (P_V^\perp)\tr,\quad \forall Z \in \Rmn,
\end{equation}
where we introduced the operators $P_U:=UU^\top M$ and $P_V:=VV^\top \Xi$ that project orthogonally onto the subspaces spanned by the columns of $U$ and $V$, respectively. With this notation, for any $H\in T_Y\cMr$, the unique triplet $(L,U_p,V_p)$ can be recovered through the relations
\begin{equation}\label{eq:tangent_parametrization}
L=U\tr M H \Xi V,\quad U_p=P_U^\perp H \Xi V,\quad V_p=P_V^\perp H\tr M U.
\end{equation}

The map $Y\mapsto \cP_Y$ is smooth, $\cMr$ being a smooth embedded submanifold of a Euclidean space (see \cite[Chapter 3]{boumal_2023}), and we denote its differential along a tangent direction $H\in T_{Y}\cMr$ by $D \cP_Y[H]\in\cL(\mathbb{R}^{m\times n})$. Since $\cP_Y$ is self-adjoint for every $Y\in \cMr$, so is $D \cP_Y[H]$.
Moreover, differentiating the identity $\cP^2_Y=\cP_Y$ along $H$ gives
\begin{equation*}\label{eq:properties_differential_projector}
(I-\cP_Y)D\cP_Y[H]=D\cP_Y[H]\cP_Y,
\end{equation*} 
which implies that 
\begin{equation}\label{eq:properties_differential_projector2}
\begin{aligned}
&D\cP_Y[H]Z\in N_{Y}\cMr,\quad \forall\; H,\; Z\in T_{Y}\cMr,\\
&D\cP_Y[H]W\in T_{Y}\cMr,\quad\forall\; H\in T_{Y}\cMr,\; \forall\; W\in N_Y\cMr.
\end{aligned} 
\end{equation}

Next, we wish to derive an explicit expression for $D\cP_Y[H]Z$, $Z\in\mathbb{R}^{m\times n}$, which measures the variation of $\cP_Y Z$ if the reference point $Y$ moves infinitesimally along the direction given by $H$. To this end, we deduce from \eqref{eq:proj-formula2} that
\begin{equation}\label{eq:derivativeProjector}
D \cP_Y[H] Z = DP_U[H] Z (P_V^\perp)\tr + P_U^\perp Z (DP_V[H])\tr.
\end{equation}
To express $DP_U[H]$ and $DP_V[H]$, we consider a curve $(-\varepsilon,\varepsilon)\ni t\mapsto Y(t)=U(t)S(t)V(t)\tr \in \cMr$ such that $Y(0)=Y$ and
\begin{equation}\label{eq:curve_on_Mr}
\dt Y|_{t=0}=\dt U S V\tr + U\dt S V\tr + US\dt V\tr= H,
\end{equation}
with the gauge conditions $\dt U\tr M U=0$ and $\dt V\tr \Xi V=0$, which can be imposed without loss of generality. Along such a curve, we have that
\begin{equation}\label{eq:directional_derivatives_PU_PV}
DP_U[H]=\dt U U\tr M + U \dt U\tr M\quad\text{and}\quad DP_V[H]=\dt V V\tr \Xi + V \dt V\tr \Xi.
\end{equation}
The expressions for $\dt U$ and $\dt V$ are 
\begin{equation}\label{eq:dotU_dotV}
\dt U = P_U^\perp H\Xi VS^{-1},\quad \dt V = P_V^\perp H\tr M US^{-\top},
\end{equation}
see \cite[Section 2]{kochlubich2007}. Inserting \eqref{eq:dotU_dotV} into \eqref{eq:directional_derivatives_PU_PV} and \eqref{eq:derivativeProjector}, we derive 
\begin{equation}\label{eq:derivative_projection}
\begin{aligned}
D \cP_Y[H] Z&\!=\!\left[ P_U^\perp H \Xi V S^{-1}U\tr + U S^{-\top} V\tr \Xi H\tr (P_U^\perp)\tr \right]M Z (P_V^\perp)\tr\\
&+P_U^\perp Z \Xi \left[P_V^\perp H\tr M U S^{-\top} V\tr + V S^{-1} U\tr M H (P_V^\perp)\tr  \right].
\end{aligned}
\end{equation}
which also admits the compact formulation
\begin{equation*}\label{eq:derivative_projection_compact}
\begin{aligned}
D \cP_Y[H] Z&=\left[ U_p S^{-1}U\tr + U S^{-\top} U_p\tr \right]M Z (P_V^\perp)\tr +P_U^\perp Z \Xi \left[V_p S^{-\top} V\tr + V S^{-1} V_p\tr  \right],
\end{aligned}
\end{equation*}
using the triplet $(L,U_p,V_p)$ that identifies uniquely the tangent vector $H$. 
Whenever $Z\in T_{Y}\cMr$, the expression simplifies since $U_p^\top M Z (P_V^\perp)^\top=0$ and 
$(P_U^\perp) Z\Xi V_p=0$, yielding
\begin{equation}\label{eq:derivative_projection_compact_tangent}
\begin{aligned}
D \cP_Y[H] Z&= U_p S^{-1} \widetilde{V}_p^\top + \widetilde{U}_p S^{-1}V_p^\top,
\end{aligned}
\end{equation}
$Z$ being parametrized by the triplet $(\widetilde{L},\widetilde{U}_p,\widetilde{V}_p)$. Equation \eqref{eq:derivative_projection_compact_tangent} is symmetric in the tangent factors of $H$ and $Z$ showing that $D\cP_{Y}[H]Z=D\cP_{Y}[Z]H$ for every $H,Z\in T_{Y}\cMr$.

Next, we recall the DLRA framework introduced in \cite{kochlubich2007,nonnenmacher2008dynamical}. Given a matrix ODE such as \eqref{eq:state_equation}, DLRA looks for an approximated solution $Y(t)$ by prescribing that its tangent vector $\dt Y(t)$ solves instantaneously the minimization problem
\begin{equation}
  \label{eq:dlra-min}
  \min_{\dt Y(t) \,\in\, T_{Y(t)} \cMr}\,
    \bigl\| F(t, Y(t), \zb(t)) - \dt Y(t) \bigr\|_{M, \Xi}^2.
\end{equation}
Due to a simple variational argument, the minimization \eqref{eq:dlra-min} is equivalent to the Galerkin
condition
\begin{equation*}
  \label{eq:dlra-galerkin}
  \langle\, \dt Y(t) - F(t, Y(t), \zb(t)),\, \delta Y(t) \,\rangle_{M, \Xi} \;=\; 0
  \qquad \forall\, \delta Y(t) \in T_{Y(t)} \cMr,\quad \text{for a.e. } t\in (0,T]
\end{equation*}
and, equivalently again, to the projected
matrix ODE 
\begin{equation}
  \label{eq:dlra-projected}
  \dt Y(t) \;=\; \cP_{Y(t)}\, F\bigl(t, Y(t),\, \zb(t) \bigr),
  \qquad Y(0) \in \cMr, \quad \text{for a.e. } t\in (0,T].
\end{equation}
The well-posedness of \eqref{eq:dlra-projected} on the manifold $\cMr$ can be analyzed through standard ODE theory in the finite-dimensional case (see, e.g., \cite{hairer2011solving}), requiring typically that the $r$-th largest singular value of $Y$ is uniformly bounded from below over $I$ in order to guarantee the Lipschitz continuity of the projection operator appearing on the right-hand side. In the infinite-dimensional case, the topic has been the subject of dedicated works, see \cite{bachmayr2021existence} for parabolic equations, and \cite{kazashi2021existence} for random semilinear evolution equations. For the sake of this manuscript we make the following working assumption.
\begin{assumption}\label{ass:well-posedness_DLRA}
For any $\zb \in L^2(I;\Rs)$, the projected state equation \eqref{eq:dlra-projected} admits a unique solution $Y\in H^1(I;\Rmn)$ with $Y(t)\in \cMr$ for every $t\in [0,T]$.
\end{assumption}

\section{Dynamical low-rank approximation and optimality systems}
\label{sec:DLRAOPT}

In this section, we present two distinct DLRA-based strategies for efficiently approximating gradients for the optimal control problem \eqref{eq:ocp}. 
\subsection{Optimize-then-Project}\label{subsec:OtP}

The central idea of this approach is to perform a dynamical low-rank approximation of both the state and adjoint equations. Denoting with $Y(t)$ and $Q(t)$ the DLRA approximated solutions of \eqref{eq:ocp_optimality_1} and \eqref{eq:ocp_optimality_2} respectively, the resulting approximated optimality system reads
\begin{align}
  \dt Y(t) &=\; \cP_{Y(t)} F\left(t, Y(t), \zb(t) \right),&\quad \text{for a.e. }t\in I,\label{eq:opt-DLRAstate}\\
  -\dt Q(t) &=\; \cP_{Q(t)}\left(F_X^\star(t,Y(t),\zb(t))Q(t)+X_d(t)-Y(t)\right),&\quad \text{for a.e. } t\in I,\label{eq:opt-DLRAadjoint}\\
   \nu \zb(t)&=\; B^\star(t,Y(t)) Q(t),&\quad \text{for a.e. } t\in I\label{eq:opt-DLRAcontrol}
\end{align}
together with the initial condition $Y(0)=Y_{0,r}$, where $Y_{0,r}$ is a rank-$r$ approximation of the initial datum $X_0$, and the final condition $Q(T)=\varepsilon_0 Q_T$, where $\varepsilon_0>0$ is arbitrarily small and $Q_T\in \Rmn$ has rank $r$ and norm equal to one.

\begin{remark}[On the final condition for $Q$]\label{rmk:finalconditionOtP}
Since $P(T)=0$, the most natural final condition would be $Q(T)=0$. However, the zero matrix does not belong to $\cMr$ for any $r\in \mathbb{N}$ and therefore it cannot serve as a valid terminal condition for \eqref{eq:opt-DLRAadjoint}. A rigorous mathematical remedy is to initialize $Q$ at time $T$ with a rank-$r$ random perturbation of the zero matrix of norm tuned by $\varepsilon_0$. This parameter will explicitly appear in the error bound on the gradient of Theorem \ref{thm:OptProj-errorestimate}. In practical computations, a common workaround is to first integrate the (full-rank) adjoint equation backward in time for a small number $n_s$ of steps. The resulting solution is then truncated to rank $r$ providing a valid starting point in $\cMr$, before switching to a DLRA integrator (see Section \ref{sec:disc}) for the remaining time steps.

Note that the above discussion is specific to the functional \eqref{eq:ocp}. If the objective functional also included a terminal-state penalty such as $\frac{\alpha}{2}\|X(T)-X_d(T)\|^2_{M,\Xi}$, then the adjoint variable $P$ would satisfy $P(T)~=~\alpha(X_d(T)-X(T))$. In that case, a natural final condition for $Q$ would be the rank-$r$ truncation of the term $\alpha(X_d(T)-Y(T))$.
\end{remark}

We emphasize that \eqref{eq:opt-DLRAstate}--\eqref{eq:opt-DLRAadjoint} are obtained by replacing the two ambient ODEs \eqref{eq:ocp_optimality_1}--\eqref{eq:ocp_optimality_2} by their DLRA-projected versions. The two approximations are independent and consequently
the dynamics are projected onto the two different tangent spaces, $T_{Y(t)} \cMr$ and $T_{Q(t)} \cMr$. In principle, it is also possible to consider a distinct rank, say $\widetilde{r}\in \N$, for the adjoint equation. As this requires only minor notational changes in what follows, we restrict to $\widetilde{r}=r$ for the rest of this work. 

\subsubsection{Error estimate for the gradient}

The OtP approach provides only an approximation to the true gradient of $\widehat{\cJ}$ and, in general, \eqref{eq:opt-DLRAstate}--\eqref{eq:opt-DLRAcontrol} may not actually correspond to the optimality system of an easily identifiable auxiliary functional. To quantify the error introduced by this approximation, we follow the framework of \cite{kieri2016robust,kieri2019projection} and we formulate some additional assumptions. The first one is on the control operator.
\begin{assumption}\label{ass:control_operator}
The control operator $B$ is time and state independent, that is, it is a given element of $\cL(\Rs,\Rmn)$.
\end{assumption}
To ease the notation, we still denote the (time and state independent) control operator by $B$, and its adjoint by $B^\star$. Note that by Assumptions \ref{ass:vector_field} and \ref{ass:control_operator}, the Fr\'echet derivative of $F$ with respect to $X$ depends only on $t$ and $X$ and not on the control $\zb$.
The second and third assumptions are related to the modelling error, that is, we ask that both $F$ and $F_X$ almost map onto the tangent bundle of $\cMr$;
\begin{assumption}\label{ass:modelling_error_F}
Let $I_Y$ and $I_{\zb}$ be suitable neighbourhoods of the trajectories $\{Y(t)\}_{t\in I}$ and $\{\zb(t)\}_{t\in I}$. Then, there exists a positive constant $\varepsilon_Y$ such that:
\begin{equation}\label{eq:modelling_assumption_F}
\|\cP_{Y}^\perp F(t,Y,\zb)\|_{M,\Xi}\leq \varepsilon_Y,\quad\text{for all } Y\in \cMr\cap I_Y,\; \zb\in I_{\zb}\text{ and for a.e. } t\in I.
\end{equation}
\end{assumption}
\begin{assumption}\label{ass:modelling_error_Fx}
Let $I_Y$ be the neighbourhood of Assumption \ref{ass:modelling_error_F} and $I_Q$ be a suitable neighbourhood of the trajectory $\{Q(t)\}_{t\in I}$. Then, there exists a positive constant $\varepsilon_Q$ such that:
\begin{equation}\label{eq:modelling_assumption_Fx}
\|\cP_{Q}^\perp (F_X^\star(t,Y)Q+X_d-Y)\|_{M,\Xi}\leq \varepsilon_Q,\quad\text{for all } Q\in \cMr\cap I_Q,\; Y\in \cMr \cap I_Y\text{ and for a.e. } t\in I.
\end{equation}
\end{assumption}
The fourth assumption is concerned with one-sided Lipschitz constants.
\begin{assumption}\label{ass:assumption_onesided_Lipschitz}
We assume that $F$ satisfies the one-sided Lipschitz bound
\begin{equation}\label{eq:one_sided_Lipschitz_bound_F}
\langle X-Y, F(t,X,\zb)-F(t,Y,\zb)\rangle_{M,\Xi}\leq \ell_F \|X-Y\|^2_{M,\Xi},\quad \text{for all }X,Y\in \Rmn,\;\zb\in \Rs \text{ and for a.e. } t\in I.
\end{equation}
Furthermore, the Fr\'echet derivative of $F$ with respect to $X$ is Lipschitz continuous with constant $L_{F_X}$ and satisfies the one-sided Lipschitz bound
\begin{equation}\label{eq:one_sided_Lipschitz_bound_Fx}
\langle P-Q, F_X^\star(t,Y)(P-Q)\rangle_{M,\Xi}\leq \ell_{F_X} \|P-Q\|^2_{M,\Xi},\quad \text{for all }P,Q\in \Rmn,\; Y\in I_Y\cap \cMr\text{ and for a.e. } t\in I.
\end{equation}
with $I_Y$ as in Assumption \ref{ass:modelling_error_F}.
\end{assumption}
Assumptions \ref{ass:modelling_error_F} and \ref{ass:assumption_onesided_Lipschitz} are sufficient to derive an error bound on the solution of the DLRA state equation (see \cite[Section 2.2]{kieri2019projection}), that is, for every $t\in [0,T]$ it holds that
\begin{equation}\label{eq:error_state}
\|X(t)-Y(t)\|_{M,\Xi}\leq e^{\ell_F t}\|X(0)-Y_{0,r}\|_{M,\Xi}+\varepsilon_Y \int_0^t e^{\ell_F (t-s)}\rd s =: e_Y(t).
\end{equation}
We are then ready to derive an estimate for the error on the gradient introduced by \eqref{eq:opt-DLRAstate}--\eqref{eq:opt-DLRAcontrol}.
\begin{theorem}\label{thm:OptProj-errorestimate}
Let Assumptions \ref{ass:vector_field}, \ref{ass:well-posedness_DLRA}, \ref{ass:control_operator}, \ref{ass:modelling_error_F}, \ref{ass:modelling_error_Fx} and \ref{ass:assumption_onesided_Lipschitz} hold and let $g(\zb):=\nu \zb -B^\star Q$ be the approximated gradient.
Then, for every $\zb\in L^2(I;\Rs)$ and for all $t\in [0,T]$ 
\begin{equation}\label{eq:estimate_error_gradient}
\|\nabla \widehat{\cJ}(\zb)(t)-g(\zb)(t)\|_\Gamma \leq  \|B\|_{\op}\left(e^{\ell_{F_X}(T-t)}\varepsilon_0 + \int_{t}^{T} e^{\ell_{F_X}(\sigma-t)}\left[\varepsilon_Q + (1+L_{F_X}C_P) e_Y(\sigma)\right]\rd \sigma\right),
\end{equation}
where $C_P:=\sup_{t\in [0,T]} \|P(t)\|_{M,\Xi}$.
\end{theorem}
\begin{proof}
We start observing that $\nabla \widehat{\cJ}(\zb)-g(\zb)=B^\star (Q-P)$, and thus $\|\nabla \widehat{\cJ}(\zb)(t)-g(\zb)(t)\|_\Gamma\leq \|B\|_{\op}R(t)$ with $R(t):=\|P(t)-Q(t)\|_{M,\Xi}$. Note that $R(T)=\varepsilon_0$ due to the final condition imposed on $Q$.
For the sake of this proof, we introduce $\tau=T-t$ and denote with a superscript tilde all quantities expressed as a function of $\tau$, i.e., $\widetilde{R}(\tau):=R(T-\tau)$.

Leveraging the evolution equations \eqref{eq:ocp_optimality_1}--\eqref{eq:ocp_optimality_3} together with \eqref{eq:opt-DLRAstate}--\eqref{eq:opt-DLRAadjoint}, and with some algebraic manipulations we obtain for a.e. $\tau$,
\begin{equation}\label{eq:proof_optthenproject_step1}
\begin{aligned}
\frac{1}{2}\frac{d}{d\tau}\widetilde{R} ^2 &=\langle \widetilde{P}-\widetilde{Q}, \dt{\widetilde{P}}-\dt{\widetilde{Q}}\rangle_{M,\Xi}\\
&= \langle \widetilde{P}-\widetilde{Q}, F_X^\star(T-\tau,\widetilde{X})\widetilde{P} +\widetilde{X_d}-\widetilde{X}  -\cP_{\widetilde{Q}}\left(F_X^\star(T-\tau,\widetilde{Y})\widetilde{Q}+\widetilde{X_d}-\widetilde{Y}\right)\rangle_{M,\Xi}\\
&=\langle \widetilde{P}-\widetilde{Q},  F_X^\star(T-\tau,\widetilde{X})\widetilde{P} -F_X^\star(T-\tau,\widetilde{Y})\widetilde{Q} \rangle_{M,\Xi}+\langle \widetilde{P}-\widetilde{Q},\widetilde{Y}-\widetilde{X}\rangle_{M,\Xi}\\
&+\langle  \widetilde{P}-\widetilde{Q}, \cP_{\widetilde{Q}}^\perp \left(F_X^\star(T-\tau,\widetilde{Y})\widetilde{Q} +\widetilde{X_d}-\widetilde{Y}\right)\rangle_{M,\Xi}\\
&=\langle \widetilde{P}-\widetilde{Q},  F_X^\star(T-\tau,\widetilde{Y})(\widetilde{P}-\widetilde{Q})\rangle_{M,\Xi} +\langle \widetilde{P}-\widetilde{Q},(F_X^\star(T-\tau,\widetilde{X})-F_X^\star(T-\tau,\widetilde{Y}))\widetilde{P}\rangle_{M,\Xi}\\
&+\langle \widetilde{P}-\widetilde{Q}, \widetilde{Y}-\widetilde{X}\rangle_{M,\Xi}+\langle \widetilde{P}-\widetilde{Q},  \cP_{\widetilde{Q}}^\perp (F_X^\star(T-\tau,\widetilde{Y})\widetilde{Q} +\widetilde{X_d}-\widetilde{Y})\rangle_{M,\Xi}
\end{aligned}
\end{equation}
We now bound the four terms separately using, respectively, \eqref{eq:one_sided_Lipschitz_bound_Fx}, the Lipschitz continuity of $F_X$,   \eqref{eq:error_state} and \eqref{eq:modelling_assumption_Fx}:
\begin{equation}\label{eq:proof_optthenproject_step2}
\begin{aligned}
\frac{1}{2}\frac{d}{d\tau}\widetilde{R} ^2 \leq  \ell_{F_X}\|\widetilde{P}-\widetilde{Q}\|_{M,\Xi}^2 +(\varepsilon_Q+ (1+L_{F_X}C_P)\widetilde{e_Y} ) \|\widetilde{P}-\widetilde{Q}\|_{M,\Xi}.
\end{aligned}
\end{equation}
Note that $C_P<\infty$ since $P\in H^1(I;\Rmn)$ and thus in $C^0(\overline{I};\Rmn)$.
Equation \eqref{eq:proof_optthenproject_step2} implies
\[\dt{\widetilde{R}} \leq \ell_{F_X} \widetilde{R} + (\varepsilon_Q+(1+L_{F_X}C_P)\widetilde{e_Y})\quad\text{for a.e. }\tau\in (0,T],\]
and Gr\"onwall's lemma finally leads to
\[\widetilde{R}(\tau)\leq e^{\ell_{F_X}\tau}\varepsilon_0 + \int_0^\tau e^{\ell_{F_X}(\tau-s)}\left[\varepsilon_Q+(1+L_{F_X}C_P)\widetilde{e_Y}(s)\right]\rd s,\quad \forall \tau \in[0,T].\]
Transforming back to the $t$ variable via the change of variable $\sigma=T-s$ yields the claim.
\end{proof}

\begin{remark}[Consistency in the full-rank limit]
\label{rem:full-rank}
The dimension of $T_{Y}\cMr$ is $(m+n-r)r$ so that the tangent space is all of $\Rmn$ if $r=\min(m,n)$.
In such a case, we have $\varepsilon_Y=\varepsilon_Q=0$ and $Y_{0,r}=X(0)$ is an admissible initial condition on $\cMr$ provided that $X(0)$ is of full-rank. The only remaining error contribution in \eqref{eq:estimate_error_gradient} is due to the terminal condition on the adjoint variable, which can be made arbitrarily small by tuning $\varepsilon_0$.
\end{remark}

\begin{remark}[Affine vector field]
\label{rem:linear_quadratic}
If the vector field $F$ is affine in $X$, then the second term in the last step of \eqref{eq:proof_optthenproject_step1} vanishes and the error bound \eqref{eq:estimate_error_gradient} simplifies to 
\begin{equation}\label{eq:estimate_error_gradient_linear}
\|\nabla \widehat{\cJ}(\zb)(t)-g(\zb)(t)\|_\Gamma \leq  \|B\|_{\op}\left(e^{\ell_{F_X}(T-t)}\varepsilon_0 + \int_{t}^{T} e^{\ell_{F_X}(\sigma-t)}\left[\varepsilon_Q + e_Y(\sigma)\right]\rd \sigma\right).
\end{equation}
\end{remark}

\subsection{Project-then-Optimize}
\label{subsec:PtO}
The starting point of the second strategy is to consider an approximation of \eqref{eq:ocp_full}--\eqref{eq:state_equation} obtained by replacing the state equation with its DLRA approximation, namely
\begin{align}
  & \min_{Y\in H^1(I;\mathbb{R}^{m\times n}),\; \zb\in L^2(I;\mathbb{R}^s)} \mathcal{J}(Y, \zb)\label{eq:ocp_full2}\\
      &\text{subject to }\; \begin{cases}
  \dt Y(t) &=\; \cP_{Y(t)} F\bigl(t, Y(t),\, \zb(t) \bigr),\quad \text{for a.e. }t\in (0,T],\\
   Y(0) &= Y_{0,r}.
  \end{cases}\label{eq:state_equationDLRA}
\end{align}

Due to Assumption \ref{ass:well-posedness_DLRA}, we can define the solution operator $\Sct: L^2(I;\Rs)\rightarrow H^1(I;\mathbb{R}^{m\times n})$ such that $\Sct(\zb)=Y$, and rewrite \eqref{eq:ocp_full2}--\eqref{eq:state_equationDLRA} as the reduced optimal control problem,
\begin{equation}\label{eq:ocp-matrix-dlra-reduced}
\min_{\zb\in L^2(I;\Rs)} \Jct(\zb),
\end{equation}
with $\Jct(\zb):=\mathcal{J}(\Sct(\zb),\zb)$.
The next lemma proves the existence of a minimizer.

\begin{lemma}\label{lem:existence_minimizer_ocp_dlra}
Let Assumptions \ref{ass:vector_field} and \ref{ass:well-posedness_DLRA} hold.
Then, \eqref{eq:ocp-matrix-dlra-reduced} admits a minimizer.
\end{lemma}
\begin{proof}
The proof follows that of Lemma \ref{lem:existence_minimizer_ocp} with some minor modifications. Let $\{\zb_k\}_{k\in\N}$ be a minimizing bounded sequence converging weakly to $\zb^\star$ and $\{Y_k\}_{k\in\N}$, with $Y_k:=\Sct(\zb_k)$, be the sequence of solutions to \eqref{eq:state_equationDLRA} with force term $\zb_k$.
Since $\|\cP_{Y(t)} F(t,Y(t),\zb(t))\|_{M,\Xi}\leq \|F(t,Y(t),\zb(t))\|_{M,\Xi}$, we can show, with the same arguments of Lemma \ref{lem:existence_minimizer_ocp}, that the sequence $\{Y_k\}_{k\in\mathbb{N}}$ is uniformly bounded and equicontinuous,
implying the existence of a subsequence that converges strongly to $Y^\star$ in $C^0(\overline{I};\mathbb{R}^{m\times n})$.
Since every $Y_k(t)$ has rank $r$, the limit $Y^\star(t)$ has rank at most $r$, and $\cP_{Y^\star(t)}$ is defined only if $Y^\star(t)\in\cMr$. We first assume that this holds for every $t\in\overline{I}$; the other case is excluded at the end of the proof.
To identify $Y^\star$ with $\Sct(\zb^\star)$, we consider the integral representation
\[Y_k(t)=Y_{0,r} +\int_0^t \cP_{Y_k(s)}\left(f(s,Y_k(s))+B(s,Y_k(s))\zb_k(s)\right) \rd s.\]
Since $Y_k$ converges strongly to $Y^\star$ in $C^0(\overline{I};\mathbb{R}^{m\times n})$ and the mapping $\cMr\ni Y\mapsto \cP_Y$ is continuously differentiable, $\cP_{Y_k(s)}f(s,Y_k(s))$ converges to $\cP_{Y^\star(s)}f(s,Y^\star(s))$ for a.e. $s$ as $k\rightarrow \infty$. Furthermore, due to the boundedness of the sequence $\{Y_k\}_{k\in \mathbb{N}}$ and the growth assumption on $f$, the dominated convergence theorem implies that
\[\lim_{k\rightarrow \infty} \int_0^t \cP_{Y_k(s)}f(s,Y_k(s))\rd s=\int_0^t \cP_{Y^\star(s)}f(s,Y^\star(s))\rd s.\]
Concerning the last term, we write
\begin{equation*}
\begin{aligned}
&\left\|\int_0^t \left(\cP_{Y_k(s)}B(s,Y_k(s))\zb_k(s)-\cP_{Y^\star(s)}B(s,Y^\star(s)) \zb^\star(s)\right)\rd s\right\|_{M,\Xi}\\ & \leq \left\|\int_0^t \left(\cP_{Y_k(s)}B(s,Y_k(s))-\cP_{Y^\star(s)}B(s,Y^\star(s))\right) \zb_k(s)\rd s\right\|_{M,\Xi}
+\left\|\int_0^t \cP_{Y^\star(s)}B(s,Y^\star(s))\left(\zb_k(s)-\zb^\star(s) \right)\rd s\right\|_{M,\Xi},
\end{aligned}
\end{equation*}
and observe that the first term tends to zero due to the uniform convergence of $Y_k$, the continuity of $Y\mapsto \cP_Y$ and the Lipschitz continuity of $B$, while the second term tends to zero due to the weak convergence of $\zb_k$.
Hence, $Y^\star$ solves \eqref{eq:dlra-projected} with control $\zb^\star$, and $Y^\star=\Sct(\zb^\star)$ by the uniqueness in Assumption \ref{ass:well-posedness_DLRA}.
The proof then follows that of Lemma \ref{lem:existence_minimizer_ocp}.

It remains to exclude that $Y^\star$ leaves $\cMr$. Suppose, by contradiction, that it does, and let $\tau\in(0,T]$ be the first time at which $Y^\star(\tau)\notin\cMr$. It exists and is positive, since $Y^\star$ is continuous, $Y^\star(0)=Y_{0,r}\in\cMr$, and the matrices of rank less than $r$ form a closed set. For $t\in[0,\tau)$, the argument above only involves $Y^\star(s)$ with $s\in[0,t]$, where $Y^\star(s)\in\cMr$, and it shows that $Y^\star$ solves \eqref{eq:dlra-projected} with control $\zb^\star$ on $[0,\tau)$, with the same initial value as $\Sct(\zb^\star)$. Two such solutions coincide as long as they stay in $\cMr$: the right-hand side of \eqref{eq:dlra-projected} is locally Lipschitz continuous in $Y$ on $\cMr$, by the smoothness of $Y\mapsto\cP_Y$ and Assumption \ref{ass:vector_field}, and Gr\"onwall's lemma applies. Hence, $Y^\star=\Sct(\zb^\star)$ on $[0,\tau)$. Next, let $t_n\in[0,\tau)$ with $t_n\rightarrow\tau$. Since $t\mapsto Y^\star(t)$ and $t\mapsto\Sct(\zb^\star)(t)$ are continuous on $\overline{I}$,
\[Y^\star(\tau)=\lim_{n\rightarrow\infty}Y^\star(t_n)=\lim_{n\rightarrow\infty}\Sct(\zb^\star)(t_n)=\Sct(\zb^\star)(\tau),\]
which is in $\cMr$ by Assumption \ref{ass:well-posedness_DLRA}. This contradicts $Y^\star(\tau)\notin\cMr$.
\end{proof}

\begin{remark}\label{rmk:nonconvexity}
Problem \eqref{eq:ocp-matrix-dlra-reduced} may not have a unique minimizer even if the vector field $F$ is linear in both $Y$ and $\zb$, in which case the original optimal control problem \eqref{eq:ocp} is linear-quadratic. Indeed, the projector $\cP_{Y}$ makes the solution map $\widetilde{\mathcal{S}}$ nonlinear. Hence, the map $\zb\mapsto \Jct(\zb)$ is generally not convex.
\end{remark}

\subsubsection{Derivation of optimality conditions}
We now derive the first-order optimality conditions satisfied by the local minima of \eqref{eq:ocp-matrix-dlra-reduced}. The gradient of $\Jct$ at a point $\zb$ satisfies
\begin{equation}\label{eq:gradient-dlra}
\langle\nabla \Jct(\zb),\vb\rangle_{L^2(I;\Rs)}=\langle \nabla\cJ(\Sct(\zb),\zb),\vb\rangle_{L^2(I;\Rs)}=\langle \nu \zb,\vb\rangle_{L^2(I;\Rs)} + \langle \Sct(\zb)-X_d, \Sct'(\zb)[\vb]\rangle_{L^2(I;\mathbb{R}^{m\times n})},
\end{equation}
where $\Sct'(\zb):L^2(I;\Rs)\rightarrow L^2(I;\mathbb{R}^{m\times n})$ denotes the Fr\'echet derivative of $\Sct$ at $\zb$ 
along the direction $\vb$, namely
\begin{equation}\label{eq:derivative_of_S}
\Sct'(\zb)[\vb]:=\frac{d}{d\varepsilon} Y(t;\zb+\varepsilon \vb)\left|_{\varepsilon=0}\right. .
\end{equation}
More concretely, $\Sct'(\zb)$ maps a direction $\vb$
to the solution of the matrix ODE
\begin{equation}\label{eq:deltaY}
\begin{cases}
\dt H(t)&= D\cP_{Y(t)}[H (t)]F(t,Y(t),\zb(t)) +\cP_{Y(t)}\left(F_X(t,Y(t),\zb(t))H(t)+F_{\zb}(t,Y(t),\zb(t))\vb(t)\right),\quad\text{for a.e. }t\in (0,T]\\
H(0)&= 0.
\end{cases}
\end{equation}
From the definition \eqref{eq:derivative_of_S}, it follows that $H(t)$ belongs to $T_{Y(t)} \cMr$, since $H(t)$ coincides with the velocity of the curve $c:(-\varepsilon',\varepsilon')\rightarrow \cMr$, $c(\tau):=Y(t;\zb+\tau \vb)$. This is also shown rigorously in the next proposition.
\begin{proposition}
The time-dependent matrix $H(t)$, solution of \eqref{eq:deltaY}, lies in the tangent space $T_{Y(t)}\cMr$ for any $t\in [0,T]$.
\end{proposition}
\begin{proof}
First, we extend the definition of the derivative of the projector operator to any ambient matrix direction $Z\in\Rmn$ by setting $\widetilde{D\cP}_Y[Z]=D\cP_Y[\cP_Y Z]$. Let now $N(t)=(I-\cP_{Y(t)})H(t)$ denote the normal component of $H(t)$. Our aim is to show that $N(t)=0$ for any $t\in [0,T]$. To this end, observe that the time derivative of $N(t)$ reads $\dt N(t) = -D\cP_{Y(t)}[\dt Y(t)]H(t)+(I-\cP_{Y(t)})\dt H(t)$.
Substituting the expression of $\dt H(t)$ and decomposing $H(t)=N(t)+\widehat{H}(t)$, $\widehat{H}(t):=\cP_{Y(t)} H(t)\in T_{Y(t)}\cMr$, yields
\begin{align*}
\dt N(t) &= -D\cP_{Y(t)}[\dt Y(t)](N(t)+\widehat{H}(t))+(I-\cP_{Y(t)})\widetilde{D\cP}_{Y(t)}[N(t)+\widehat{H}(t)]F(t,Y,\zb)\\
&=-D\cP_{Y(t)}[\dt Y(t)](N(t)+\widehat{H}(t))+(I-\cP_{Y(t)})D\cP_{Y(t)}[\widehat{H}(t)]F(t,Y,\zb),
\end{align*}
for a.e. $t\in [0,T]$.
Next, we use the linearity of $H\mapsto D\cP_{Y(t)}[H]$ and the symmetry of $D\cP_{Y(t)}$ on tangent vectors (see \eqref{eq:derivative_projection_compact_tangent}) to write
\begin{equation*}
\begin{aligned}
\dt N(t)&=-D\cP_{Y(t)}[\dt Y(t)]N(t)-D\cP_{Y(t)}[\widehat{H}(t)]\dt Y(t)+(I-\cP_{Y(t)})D\cP_{Y(t)}[\widehat{H}(t)]F(t,Y,\zb)
\end{aligned}
\end{equation*}
for a.e. $t\in [0,T]$.
From \eqref{eq:properties_differential_projector2}, it follows that 
\[(I-\cP_{Y(t)})D\cP_{Y(t)}[\widehat{H}(t)]F(t,Y,\zb)=(I-\cP_{Y(t)})D\cP_{Y(t)}[\widehat{H}(t)]\dt Y(t)=D\cP_{Y(t)}[\widehat{H}(t)]\dt Y(t),\] 
which inserted back leads to
\[\dt N(t)=-D\cP_{Y(t)}[\dt Y(t)]N(t),\quad \text{for a.e. }t\in [0,T],\]
which is a linear homogeneous matrix ODE for $N(t)$. Since $N(0)=0$, the uniqueness of its solution implies that $N(t)=0$ for every $t\in [0,T]$ as claimed.
\end{proof}
We now explicitly compute the gradient of $\Jct$. To this end, we calculate the expression of the adjoint operator of $\Sct'(\zb)$, and let $P\in H^1(I;\mathbb{R}^{m\times n})$ be arbitrary. For every $W\in L^2(I;\mathbb{R}^{m\times n})$ we consider
\begin{equation}
\begin{aligned}
\langle W, \Sct'(\zb)[\vb]&\rangle_{L^2(I;\mathbb{R}^{m\times n})}=\langle W, H\rangle_{L^2(I;\mathbb{R}^{m\times n})}\\
&+\langle P, \dt H - D\cP_{Y}\left[H\right]F(\cdot,Y,\zb)-\cP_{Y}\left(F_X(\cdot,Y,\zb) H+ F_{\zb}(\cdot,Y,\zb)\vb\right)\rangle_{L^2(I;\mathbb{R}^{m\times n})}.
\end{aligned}
\end{equation}
By integrating by parts with respect to time the term involving $\dt H$ and using the self-adjointness of $\cP_{Y}$, we obtain
\begin{equation*}
\begin{aligned}
\langle W, \Sct'(\zb)[\vb]\rangle_{L^2(I;\mathbb{R}^{m\times n})}&=
\langle W, H \rangle_{L^2(I;\mathbb{R}^{m\times n})}
-\langle \dt P + F_X(\cdot,Y,\zb)^\star \cP_{Y} P, H\rangle_{L^2(I;\mathbb{R}^{m\times n})} \\
&-\langle P, D\cP_{Y}[H]F(\cdot,Y,\zb)\rangle_{L^2(I;\mathbb{R}^{m\times n})}- \langle F_{\zb}(\cdot,Y,\zb)^\star \cP_Y P ,\vb\rangle_{L^2(I;\Rs)}+\langle P(T),H(T)\rangle_{M,\Xi}.
\end{aligned}
\end{equation*}
Furthermore, since $H(t)$ belongs to $T_{Y(t)}\cMr$, we can equivalently write
\begin{equation}\label{eq:intermediate_step_adjoint}
\begin{aligned}
\langle W, \Sct'(\zb)[\vb]\rangle_{L^2(I;\mathbb{R}^{m\times n})}&=
\langle \cP_{Y}W, H \rangle_{L^2(I;\mathbb{R}^{m\times n})}
-\langle \cP_{Y}\dt P+ \cP_{Y} F_X(\cdot,Y,\zb)^\star \cP_{Y} P , H\rangle_{L^2(I;\mathbb{R}^{m\times n})} \\
&-\langle P, D\cP_{Y}[H]F\rangle_{L^2(I;\mathbb{R}^{m\times n})}- \langle F_{\zb}(\cdot,Y,\zb)^\star \cP_Y P ,\vb\rangle_{L^2(I;\Rs)}+\langle \cP_{Y(T)} P(T),H(T)\rangle_{M,\Xi}.
\end{aligned}
\end{equation}
The duality pairing with $\vb$ in \eqref{eq:intermediate_step_adjoint} involves only the projection of $P$ onto the tangent space. Therefore, we would like to derive a closed equation for $Q(t):=\cP_{Y(t)} P(t)$. To this end, we decompose $F=\dt Y+ F^\perp$, with $F^\perp:=(I-\cP_Y)F\in N_{Y}\cMr$. Then,
\begin{equation}\label{eq:intermediate_step_adjoint2}
\begin{aligned}
\langle D\cP_{Y}[H] F, P\rangle_{L^2(I;\mathbb{R}^{m\times n})}&=\langle D\cP_{Y}[H](\dt Y+F^\perp), P\rangle_{L^2(I;\mathbb{R}^{m\times n})}=\langle D\cP_{Y}[\dt Y]H + D\cP_{Y}[H]F^\perp, P\rangle_{L^2(I;\mathbb{R}^{m\times n})}\\
&=\langle H, D\cP_{Y}[\dt Y] P\rangle_{L^2(I;\Rmn)} + \langle D\cP_{Y}[H]F^\perp, Q\rangle_{L^2(I;\mathbb{R}^{m\times n})},
\end{aligned}
\end{equation}
where we used the symmetry of $D\cP_{Y}$, the self-adjointness of $D\cP_{Y}[\dt Y]$ and \eqref{eq:properties_differential_projector2}.
Inserting \eqref{eq:intermediate_step_adjoint2} into \eqref{eq:intermediate_step_adjoint} and using the relation
\[\dt Q=\cP_{Y}\dt P+ D\cP_{Y}[\dt Y]P,\]
yields 
\begin{equation}\label{eq:intermediate_step_adjoint4}
\begin{aligned}
\langle W, \Sct'(\zb)[\vb]\rangle_{L^2(I;\mathbb{R}^{m\times n})}&=\langle \cP_{Y}W, H \rangle_{L^2(I;\mathbb{R}^{m\times n})}+\langle \cP_{Y(T)} P(T),H(T)\rangle_{M,\Xi}-\langle F_{\zb}(\cdot, Y,\zb)^\star Q,\vb\rangle_{L^2(I;\mathbb{R}^m)}\\
&-\langle \cP_Y \dt Q+ \cP_{Y} F_X(\cdot,Y,\zb)^\star Q +D\cP^\star_{Y,F^\perp(\cdot,Y,\zb)}Q, H\rangle_{L^2(I;\mathbb{R}^{m\times n})},
\end{aligned}
\end{equation}
where $D\cP^\star_{Y,F^\perp}$ is the adjoint operator such that $\langle D\cP_{Y}[H]F^\perp,Q\rangle_{M,\Xi}=\langle H,D\cP^\star_{Y,F^\perp} Q\rangle_{M,\Xi}$ for any $H,Q\in T_{Y}\cMr$. Its expression is given by 
\begin{equation}\label{eq:expression_adjointderivative}
D\cP^\star_{Y,F^\perp}Q = F^\perp \Xi \widetilde{V}_p S^{-\top}V\tr + US^{-\top} \widetilde{U}_p\tr M F^\perp,
\end{equation} 
$Q$ being parametrized by $(\widetilde{L},\widetilde{U}_p,\widetilde{V}_p)$, and satisfies $D\cP^\star_{Y,F^\perp} Q\in T_{Y}\cMr$,
see Appendix \ref{appendix:adjoint} for an explicit derivation. 
We now choose $Q(t)\in T_{Y(t)}\cMr$ as the solution of
\begin{equation}\label{eq:adjoint_equation_first}
\begin{cases}
\cP_Y \dt Q + \cP_{Y} F_X(\cdot,Y,\zb)^\star Q +D\cP^\star_{Y,F^\perp(\cdot,Y,\zb)}Q= \cP_{Y} W,\\
Q(T)=0.
\end{cases}
\end{equation}
Substituting this choice of $Q$ into the previous identity yields
\begin{align*}
\langle W, \Sct'(\zb)[\vb]\rangle_{L^2(I;\mathbb{R}^{m\times n})}&=-\langle F_{\zb}(\cdot,Y,\zb)^\star Q, \vb\rangle_{L^2(I;\Rs)},
\end{align*}
which inserted into \eqref{eq:gradient-dlra} gives the expression of $\nabla \widetilde{\cJ}(\zb)$.
We summarize the above calculations in the following proposition.
\begin{proposition}
Let Assumptions \ref{ass:vector_field} and \ref{ass:well-posedness_DLRA} hold. Furthermore, assume that the vector field $F$ is Fr\'echet differentiable with respect to $X$ and $\zb$. Then, the stationary points $(Y,Q,\zb)\in H^1(I;\Rmn)\times H^1(I;\Rmn)\times L^2(I;\Rs)$ of \eqref{eq:ocp_full2}--\eqref{eq:state_equationDLRA} satisfy the optimality system,
\begin{align}
\medmuskip=-1mu
\thinmuskip=-1mu
\thickmuskip=-1mu
\nulldelimiterspace=0.9pt
\scriptspace=0.9pt 
\arraycolsep0.9em 
\dt Y(t) &= \cP_{Y(t)} F\bigl(t, Y(t),\, \zb(t) \bigr),&\; \text{for a.e. }t\in I,\label{eq:opt-DLRAPstate}\\
-\cP_{Y(t)} \dt Q(t) &= \cP_{Y(t)}\left(F_X(t,Y(t),\zb(t))^\star Q(t)+X_d(t)-Y(t)\right)+ D\cP^\star_{Y(t),F^\perp(t,Y(t),\zb(t))}Q(t),&\; \text{for a.e. }t\in I,\label{eq:opt-DLRAPadjoint}\\
\nu \zb(t) &= F_{\zb}(t,Y(t),\zb(t))^\star Q(t),&\; \text{for a.e. }t\in I,\label{eq:opt-DLRAPcontrol}
\end{align}
together with $Q(t)\in T_{Y(t)}\cMr$ for every $t\in [0,T)$ and the initial and final conditions $Y(0)=Y_{0,r}$ and $Q(T)=0$.
\end{proposition}

Compared to \eqref{eq:opt-DLRAstate}--\eqref{eq:opt-DLRAcontrol}, only the adjoint equation changes and the adjoint variable $Q(t)$ plays a different role: it is not a point on $\cMr$, but a tangent vector attached to the moving point $Y(t)$. 
As $Y(t)$ moves, the tangent space $T_{Y(t)}\cMr$ varies with time. Consequently, \eqref{eq:opt-DLRAPadjoint} is naturally interpreted in terms of the covariant derivative (see, e.g., \cite[Chapter 15]{gallier2020differential}), which describes the evolution of $Q(t)$ as observed in a frame moving with the tangent space $T_{Y(t)}\cMr$. Indeed, as all terms on the right-hand side of \eqref{eq:opt-DLRAPadjoint} live in $T_{Y(t)}\cMr$, \eqref{eq:opt-DLRAPadjoint} determines only the tangential part of $\dt Q(t)$ and must therefore be complemented by the explicit constraint $Q(t)\in T_{Y(t)}\cMr$ for every $t\in [0,T)$. This is a \textit{covariant} formulation.

Note that, since the tangent plane rotates, $\dt Q(t)$ does not in general lie in $T_{Y(t)}\cMr$. This is easily verified by differentiating in time the relation $Q(t)=\cP_{Y(t)}Q(t)$ which yields
\begin{equation}\label{eq:derivativetimeQ}
\dt Q(t) = \cP_{Y(t)}\dt Q(t) + D\cP_{Y(t)}[\dt Y(t)]Q(t),
\end{equation}
with $D\cP_{Y(t)}[\dt Y(t)]Q(t)\in N_{Y(t)}\cMr$ due to \eqref{eq:properties_differential_projector2}.
An equivalent \textit{ambient} formulation follows from \eqref{eq:derivativetimeQ} and reads
\begin{equation}\label{eq:opt-DLRAPadjointambient}
-\dt Q(t) = \cP_{Y(t)}\left(F_X(t,Y(t),\zb(t))^\star Q(t)+X_d(t)-Y(t)\right)+ D\cP^\star_{Y(t),F^\perp(t,Y(t),\zb(t))}Q(t) -D\cP_{Y(t)}[\dt Y(t)]Q(t),
\end{equation}
for a.e. $t\in I$. The two formulations are mathematically equivalent, but \eqref{eq:opt-DLRAPadjoint} is suitable for a DLRA-type integration, hopefully yielding computational savings, while \eqref{eq:opt-DLRAPadjointambient} needs to be integrated with any standard matrix ODE solver. The ambient formulation will nevertheless be useful for the theoretical analysis.

\subsubsection{Error estimate for the gradient}

We now provide error estimates for the approximation introduced by \eqref{eq:opt-DLRAPstate}--\eqref{eq:opt-DLRAPcontrol} compared to \eqref{eq:ocp_optimality_1}--\eqref{eq:ocp_optimality_3}. To do so, we start by estimating the second term in the right-hand side of \eqref{eq:opt-DLRAPadjoint}.

\begin{lemma}
\label{lem:curvature_normal}
Let Assumption \ref{ass:modelling_error_F} hold, $Y=USV\tr\in\cMr$, $Q\in T_Y\cMr$ be parametrized by the triplet $(L,U_p,V_p)$, and let $F:=F(t,Y,\zb)$ for some $t\in I$, $\zb\in I_{\zb}$ and $Y\in\cMr\cap I_Y$, with
$F^\perp:=(I-\cP_Y)F$. Furthermore, denote by $\sigma_r(Y)$ the smallest singular value of $S$. Then, the following bound holds
\begin{align}
  \bigl\|D\cP^\star_{Y,F^\perp}\,Q\bigr\|_{M,\Xi}
  &\leq \frac{\varepsilon_Y \|Q\|_{M,\Xi}}
               {\sigma_r(Y)} .\label{eq:adjoint_curvature_bound}
\end{align}
\end{lemma}
\begin{proof}
Recall from \eqref{eq:expression_adjointderivative} that
\[D\cP^\star_{Y,F^\perp}Q = \underbrace{F^\perp \Xi V_p S^{-\top}}_{:=\widehat{U}_p}V\tr + U\underbrace{S^{-\top} U_p\tr M F^\perp}_{=:\widehat{V}_p\tr},\]
where we highlighted the factors of the tangent parametrization. Due to the gauge conditions on $\widehat{U}_p$ and $\widehat{V}_p$, it holds that
\[\|D\cP^\star_{Y,F^\perp}Q\|^2_{M,\Xi}=\|F^\perp \Xi V_p S^{-\top}\|^2_M+\|(F^\perp)^\top M U_p S^{-1}\|^2_{\Xi},\]
and the properties of the Frobenius norm lead to
\begin{equation}
\begin{aligned}
\|D\cP^\star_{Y,F^\perp}Q\|^2_{M,\Xi}&=\|(M^{\frac{1}{2}}F^\perp \Xi^{\frac{1}{2}})(\Xi^{\frac{1}{2}}V_p) S^{-\top}\|^2_F+\|(\Xi^{\frac{1}{2}}(F^\perp)\tr M^{\frac{1}{2}})(M^{\frac{1}{2}}U_p) S^{-1}\|^2_F\\
&\leq \|M^{\frac{1}{2}}F^\perp \Xi^{\frac{1}{2}}\|^2_F\|\Xi^{\frac{1}{2}}V_p\|^2_F\|S^{-\top}\|^2_2+\|\Xi^{\frac{1}{2}}(F^\perp)\tr M^{\frac{1}{2}}\|^2_F\|M^{\frac{1}{2}}U_p\|^2_F\|S^{-1}\|^2_2\\
&\leq \frac{1}{\sigma^2_r(Y)}\|F^\perp\|^2_{M,\Xi}\|Q\|^2_{M,\Xi}\leq \frac{\varepsilon_Y^2}{\sigma^2_r(Y)}\|Q\|^2_{M,\Xi}.
\end{aligned}
\end{equation}
Taking the square root yields the claim. 
\end{proof}

The next theorem provides an estimate for the error associated to the gradient of \eqref{eq:ocp-matrix-dlra-reduced}. As the projection of $F_X$ is now performed onto the tangent space at $Y(t)$, we need a new assumption which replaces Assumption \ref{ass:modelling_error_Fx}.
\begin{assumption}\label{ass:modelling_error_Fx2}
Let $I_Y$ be the neighbourhood of Assumption \ref{ass:modelling_error_F} and let $I_Q$ be a neighbourhood of the trajectory $\{Q(t)\}_{t\in I}$. There exists a positive constant $\widehat{\varepsilon}_Q$ such that:
\begin{equation}\label{eq:modelling_assumption_Fx2}
\|\cP_{Y}^\perp (F_X^\star(t,Y)Q+X_d-Y)\|_{M,\Xi}\leq \widehat{\varepsilon}_Q,\quad\forall \;Y\in \cMr \cap I_Y,\;Q\in T_Y\cMr\cap I_Q,\text{ and for a.e. } t\in I.
\end{equation}
\end{assumption}

\begin{theorem}\label{thm:ProjthenOpt-errorestimate}
Let Assumptions \ref{ass:vector_field}, \ref{ass:well-posedness_DLRA}, \ref{ass:control_operator}, \ref{ass:modelling_error_F}, \ref{ass:assumption_onesided_Lipschitz} and \ref{ass:modelling_error_Fx2} hold. Furthermore, denote by $\sigma_r(Y(t))$ the smallest singular value of the factor $S(t)$ of $Y(t)$, $t\in I$. Then, for every $\zb\in L^2(I;\Rs)$ and for all $t\in [0,T]$
\begin{equation}\label{eq:estimate_error_gradient_retract-opt}
\|\nabla \widehat{\cJ}(\zb)(t)-\nabla \widetilde{\cJ}(\zb)(t)\|_\Gamma \leq  \|B\|_{\op}\int_{t}^T e^{\ell_{F_X}(\sigma-t)}\left[\widehat{\varepsilon}_Q+(1+L_{F_X}C_P)e_Y(\sigma)+C_Q\left(\kappaeff(\sigma) + \frac{\varepsilon_Y}{\sigma_r(Y(\sigma))}\right)\right]\rd \sigma,
\end{equation}
where $\kappaeff(t):=\|D\cP_{Y(t)}[\dt Y(t)]\|_{\op}$, $C_P:=\sup_{t\in [0,T]}\|P(t)\|_{M,\Xi}$ and $C_Q:=\sup_{t\in [0,T]} \|Q(t)\|_{M,\Xi}$.
\end{theorem}
\begin{proof}
As in the proof of Theorem \ref{thm:OptProj-errorestimate}, we define $R(t):=\|P(t)-Q(t)\|_{M,\Xi}$, set $\tau=T-t$ and denote with a superscript tilde all quantities expressed as functions of $\tau$. Using \eqref{eq:ocp_optimality_1}--\eqref{eq:ocp_optimality_3} together with \eqref{eq:opt-DLRAPstate}, \eqref{eq:opt-DLRAPcontrol} and the ambient adjoint equation \eqref{eq:opt-DLRAPadjointambient}, we obtain for a.e. $\tau$
\begin{equation*}
\begin{aligned}
\frac{1}{2}\frac{d}{d\tau}\widetilde{R} ^2 &=\langle \widetilde{P}-\widetilde{Q}, \dt{\widetilde{P}}-\dt{\widetilde{Q}}\rangle_{M,\Xi}\\
&= \langle \widetilde{P}-\widetilde{Q}, F_X^\star(T-\tau,\widetilde{X})\widetilde{P} +\widetilde{X_d}-\widetilde{X}  -\cP_{\widetilde{Y}}\left(F_X^\star(T-\tau,\widetilde{Y})\widetilde{Q}+\widetilde{X_d}-\widetilde{Y}\right) \rangle_{M,\Xi}\\
&- \langle \widetilde{P}-\widetilde{Q}, D\cP^\star_{\widetilde{Y},\widetilde{F}^\perp}\widetilde{Q} -D\cP_{\widetilde{Y}}[\dt Y]\widetilde{Q} \rangle_{M,\Xi}\\
&=\langle \widetilde{P}-\widetilde{Q}, \widetilde{Y}-\widetilde{X}\rangle_{M,\Xi} +\langle \widetilde{P}-\widetilde{Q},(F_X^\star(T-\tau,\widetilde{X})-F_X^\star(T-\tau,\widetilde{Y}))\widetilde{P}\rangle_{M,\Xi}\\
&+ \langle \widetilde{P}-\widetilde{Q},  F_X^\star(T-\tau,\widetilde{Y})(\widetilde{P}-\widetilde{Q})\rangle_{M,\Xi}  +\langle \widetilde{P}-\widetilde{Q},  \cP_{\widetilde{Y}}^\perp (F_X^\star(T-\tau,\widetilde{Y})\widetilde{Q} +\widetilde{X_d}-\widetilde{Y})\rangle_{M,\Xi}\\
&-\langle \widetilde{P}-\widetilde{Q}, D\cP^\star_{\widetilde{Y},\widetilde{F}^\perp}\widetilde{Q} -D\cP_{\widetilde{Y}}[\dt Y]\widetilde{Q} \rangle_{M,\Xi}.
\end{aligned}
\end{equation*}
We now bound the five terms separately using, respectively, \eqref{eq:one_sided_Lipschitz_bound_Fx}, the Lipschitz continuity of $F_X$,   \eqref{eq:error_state}, \eqref{eq:modelling_assumption_Fx2}, and the definition of $\kappaeff$ together Lemma \ref{lem:curvature_normal},
\begin{equation}
\begin{aligned}
\frac{1}{2}\frac{d}{d\tau}\widetilde{R}^2 \leq  \ell_{F_X}\widetilde{R}^2 +\left(\widehat{\varepsilon}_Q+ (1+L_{F_X}C_P)\widetilde{e_Y} +C_Q\left(\kappaeff + \frac{\varepsilon_Y}{\sigma_r(Y)}\right)\right)  \widetilde{R}.
\end{aligned}
\end{equation}
The claim follows by Gr\"onwall's lemma as in the proof of Theorem \ref{thm:OptProj-errorestimate}.
\end{proof}

\begin{remark}[Consistency in the full-rank limit]
\label{rem:full-rank_PtO}
If $r=\min(m,n)$, then $\varepsilon_Y=\widehat{\varepsilon}_Q=0$, and the projection operator is the identity so that also $\kappaeff(t)\equiv 0$ for any $t$. Hence, the error estimate \eqref{eq:estimate_error_gradient_retract-opt} implies that $\nabla \widehat{\mathcal{J}}(\zb)=\nabla \widetilde{\mathcal{J}}(\zb)$.
\end{remark}

\begin{remark}[Dependency on the curvature]
Compared with the estimate for the OtP approach, the bound \eqref{eq:estimate_error_gradient_retract-opt} for the PtO formulation contains two additional geometric terms involving the curvature of the low-rank manifold. The first one is $\kappaeff(t)$
which measures the rate at which the tangent space $T_{Y(t)}\cMr$ rotates along the state trajectory. The second is the ratio
$\frac{\varepsilon_Y}{\sigma_r(Y(t))}$. Both terms reflect the fact that the PtO approach is not simply a projection of the dynamics, but a dynamic approximation on a given curve on the manifold with nontrivial curvature.
Note that a standard upper bound for the norm of $D\cP_{Y(t)}$ is of order $1/\sigma_r(Y(t))$ and it is sharp, in the sense that, given $Y(t)$, one can always construct tangent directions $H\in T_{Y(t)}\cMr$ for which the bound is attained \cite[Appendix A]{kieri2019projection}. The quantity $\kappaeff$ is however weaker and more relevant, as it measures the effective curvature experienced by the solution along the path $Y(t)$, rather than the maximal curvature among all possible tangent directions.
\end{remark}

\section{Time discretization}
\label{sec:disc}

We here discuss time discretization strategies for the optimality systems.
Given $N_t\in \N$, we introduce the uniform grid
$t_i = i\,\Delta t$, $I_i:=(t_i,t_{i+1}]$ for $i=0,\dots,N_t-1$, with $\Delta t := T/N_t$.  
For a current control iterate $\zb^k$, the gradient computation requires both a forward sweep that yields a discretized state solution
$\{Y_i\}_{i=0}^{N_t}$, with $Y_i\in \cMr$, and a backward sweep yielding a discretized adjoint variable $\{P_i\}_{i=0}^{N_t}$, which may either belong to $\cMr$ at each time instant (in the OtP approach) or to the tangent space at the current state approximation $Y_i$ (in the PtO approach).  

\subsection{Forward sweep}

The forward sweep is the same for both approaches, and any time integrator method suitable for DLRA equations could be used. This has been a very hot topic of research lately with many contributions among which we mention the projector-splitting (KSL) method \cite{lubichoseledets2014}, \textit{BUG} integrators \cite{ceruti2022unconventional,ceruti2022rank,nobile2026high} and projected Runge--Kutta methods \cite{kieri2019projection,carrel2024projected_exp}. In our numerical experiments of Section \ref{sec:num} we rely on the BUG integrator presented in \cite{ceruti2022unconventional}, and briefly recalled hereafter.
The BUG method
integrates \eqref{eq:dlra-projected} over the subinterval $I_i$ by first updating the basis for the spaces of columns and rows, and then by computing a new core matrix via a Galerkin projection. Specifically, let $Y_i=U_iS_iV_i\tr$ be an approximation of $Y(t_i)$. The BUG integrator performs the steps,
\begin{align}
  \textbf{K-step:}&\quad
    \text{Solve } \dt K(t) = F\bigl(t,K(t)\, V_i^\top,\zb(t)\bigr) \Xi V_i,\quad &t\in I_i,\quad K(t_i)=U_iS_i,\label{eq:K_step}\\
       &\quad \text{Set } (U_{i+1}, R_{i+1})= \mathrm{QR}_{M}(K(t_{i+1})), \\
     \textbf{L-step:}&\quad
   \text{Solve }\dt L(t) = F\bigl(t,U_{i} L(t)^\top,\zb(t)\bigr)^{\!\top} M\, U_{i},\quad &t\in I_i,\quad L(t_i)=V_iS_i\tr,\label{eq:L_step}\\
   &\quad \text{Set }(V_{i+1}, \widetilde{R}_{i+1}^\top) = \mathrm{QR}_\Xi(L(t_{i+1})),\\
    \textbf{S-step:}&\quad
  \text{Solve } \dt S(t) = U_{i+1}^\top F\bigl(t,U_{i+1} S(t)\, V_{i+1}^\top,\zb(t)\bigr)\, \Xi\, V_i,\quad &t\in I_i,\quad S(t_i)=U_{i+1}\tr U_{i} S_i V_i\tr V_{i+1}, \label{eq:S_step}
\end{align}  
and yields $Y_{i+1}= U_{i+1}S_{i+1}V_{i+1}\tr\approx Y(t_{i+1})$. Here, $\mathrm{QR}_{M}(K)$ denotes the $M$-weighted QR algorithm that factorizes $K = U R$ with $U^\top M U = I_r$. Symmetrically,
$\mathrm{QR}_\Xi(L)$ denotes the $\Xi$-weighted QR method. 
Each step is a small matrix ODE on $\mathbb{R}^{m \times r}$,
$\mathbb{R}^{n \times r}$, and $\mathbb{R}^{r \times r}$, respectively, and can be integrated with, e.g., one step of explicit
Euler or an explicit Runge--Kutta method. However, regardless of the accuracy of the inner solver, the steps \eqref{eq:K_step}--\eqref{eq:S_step} lead to a first-order in time integrator \cite[Theorem 4]{ceruti2022unconventional}.
More advanced time integrators could be used, including high-order and rank-adaptive strategies, see, e.g., \cite{nobile2026high} for a recent overview.  

\subsection{Backward sweep}
\label{sec:adj-step}

The computation of the adjoint variable in the OtP approach can be performed as the forward sweep, by using the BUG integrator \eqref{eq:K_step}--\eqref{eq:S_step} to integrate backward in time \eqref{eq:opt-DLRAadjoint}, with the additional care of addressing the terminal condition $Q(T)=0$ as described in Remark \ref{rmk:finalconditionOtP}.
In contrast, the PtO approach requires discussing the numerical integration of the covariant ODE \eqref{eq:opt-DLRAPadjoint} which is detailed in the next subsection.

\subsubsection{Time integration of a covariant ODE}\label{subsubsec:covariantODE}
\begin{figure}
\centering
\includegraphics[scale=0.45]{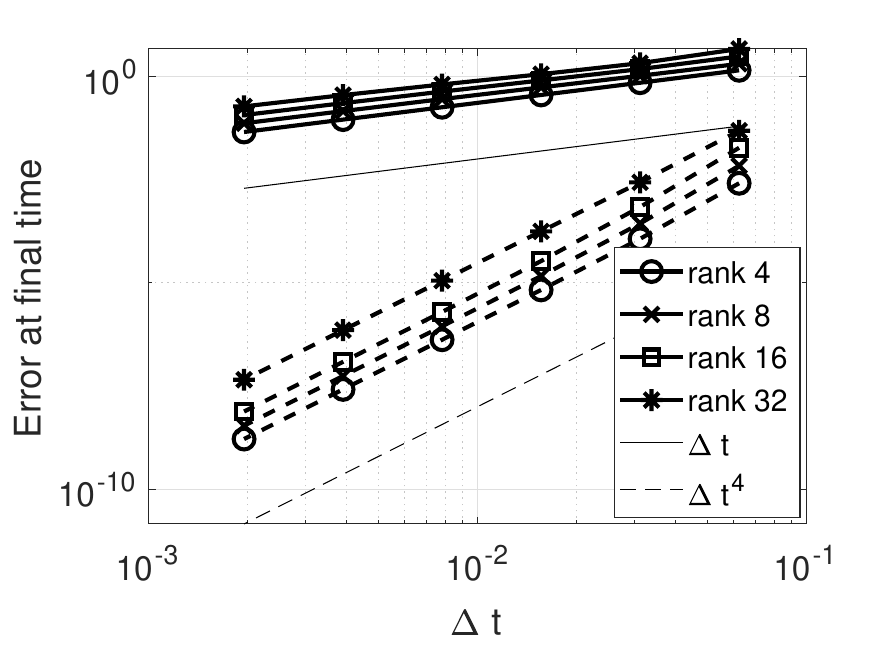}
\includegraphics[scale=0.45]{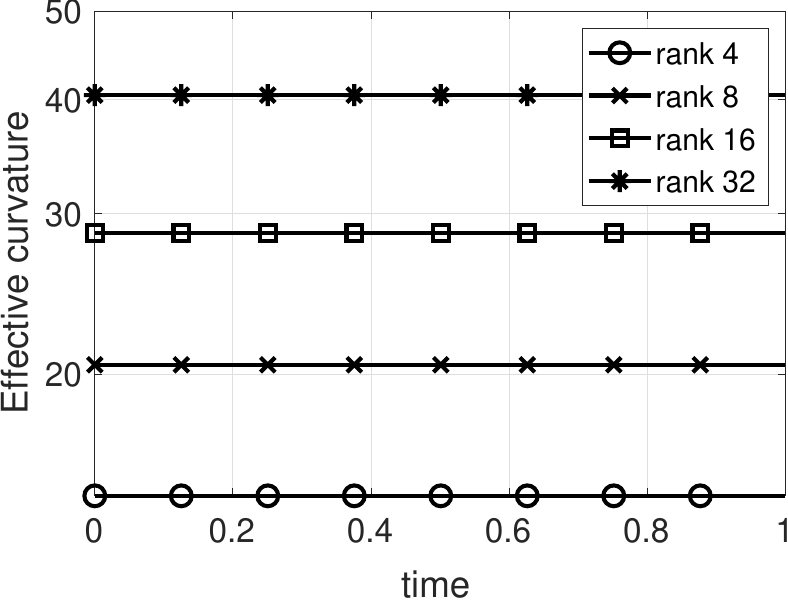}
\caption{On the left panel, decay of the error at the final time $T=1$ between the reference vector field and two approximations computed by the explicit Runge--Kutta method of order 4 applied to \eqref{eq:factorM}--\eqref{eq:factorVp} (dashed line) and the projected Euler method \eqref{eq:M_{i+1}}--\eqref{eq:U_{p,i+1}} (continuous line). On the right panel, time evolution of the effective curvature for different ranks. For this particular construction of $A_r(t)$, the effective curvature is constant in time. Parameter: $n_1=50$.}
\label{fig:NumTestCovariant}
\end{figure}

We consider the general formulation
\begin{equation}\label{eq:covariant_ODE}
\begin{cases}
&\cP_{Y(t)}\dt Q(t)=G(t,Q(t)),\quad t\in (0,T],\\
&Q(0)= Q_0,
\end{cases}
\end{equation}
where $Y$ is a curve on $\cMr$, $Q$ is a vector field along $Y$, $G$ maps to $T_{Y(t)}\cMr$, and $Q_0\in T_{Y(0)}\cMr$. A possible approach consists in expressing $Q(t)$ in its tangent parametrization 
\[Q(t)=U(t)L(t)V(t)^\top+U_p(t)V(t)^\top+U(t)V_p(t)^\top,\quad t\in (0,T],\]
and in deriving equations for the factors $L(t), U_p(t), V_p(t)$. To do so, we take the time derivative of $Q(t)$ and calculate its projection onto $T_{Y(t)}\cMr$,
\begin{equation}\label{eq:expression_projection_time_derivative}
\cP_{Y}\dt Q=U\left(\dt L+U^\top M\dt U_p+\dt V_p^\top \Xi V\right)V^\top+\left(\dt U L+(I-UU^\top M)\dt U_p\right) V^\top + U\left(\dt V L^\top + (I-VV^\top\Xi\dt) V_p\right)^\top.
\end{equation}
Here, we have additionally assumed that $Y$ is expressed via factors $U,S$ and $V$ satisfying $U^\top M\dt U=V^\top\Xi \dt V=0$. This guarantees that the tangent parametrization \eqref{eq:expression_projection_time_derivative} satisfies the gauge conditions. Next, by substituting \eqref{eq:expression_projection_time_derivative} into \eqref{eq:covariant_ODE} and imposing the equality between the factors of the tangent parametrizations, we derive the system of ODEs
\begin{align*}
\dt L+U^\top M\dt U_p+\dt V_p^\top \Xi V &= U^\top M G \Xi V,\\
\dt U L+(I-UU^\top M)\dt U_p & = (I-UU^\top M) G \Xi V,\\
\dt V L^\top + (I-VV^\top\Xi)\dt V_p &= (I-VV^\top\Xi)G^\top M U.
\end{align*}
By using the relations $U^\top M\dt U_p = -\dt U^\top M U_p$ and $\dt{V}_p^\top \Xi V= -V_p^\top \Xi \dt V$ and some algebraic manipulations, we finally derive
\begin{align}
\dt L &= U^\top M G \Xi V + \dt U^\top M U_p + V_p^\top \Xi \dt V,\label{eq:factorM}\\
\dt U_p &=(I-UU^\top M)G\Xi V -\dt U L -U\dt U^\top M U_p,\label{eq:factorUp}\\
\dt V_p &=(I-VV^\top \Xi)G^\top MU-\dt V L^\top -V\dt V^\top \Xi V_p.\label{eq:factorVp}
\end{align}
In principle, if one had an analytical expression for the factors $U$ and $V$ (as well as for their time derivatives), \eqref{eq:factorM}--\eqref{eq:factorVp} could be solved with any time integration scheme, e.g., an explicit Runge--Kutta method. We emphasize that \eqref{eq:factorM}--\eqref{eq:factorVp} \textit{do not} depend on the singular values of $Y$, but only on the effective rotation speed of the two subspaces $U$ and $V$! Consequently, we do not expect any restriction on the time step due to small singular values of $Y$ (see Figure \ref{fig:NumTestCovariant}).
Unfortunately though, in our setting we only have available (approximations of) $Y$ and its factors at prescribed time instants. One could think about replacing the time derivatives of $U$ and $V$ in \eqref{eq:factorM}--\eqref{eq:factorVp} with classical increments. This would require to first align two consecutive factors by solving a Procrustes problem (\cite[Chapter 12.4.1]{golub2013matrix}) since the BUG integrator may introduce arbitrary rotations such that $\|M^{\frac{1}{2}}(U_{i+1}-U_i)\|_{F}$ may not be of order $\mathcal{O}(\Delta t)$.

It turns out that a more promising idea is to consider the projected Euler method,
\begin{equation}\label{eq:projected_Euler}
Q_{i+1}=\cP_{Y_{i+1}}\left(Q_i+\Delta t G(t_i,Q_i)\right),
\end{equation}
which, introducing the shorthand notation $\widetilde{Q}_i=Q_i+\Delta t G(t_i,Q_i)$, leads to the following updates for the tangent factors (see \eqref{eq:tangent_parametrization})
\begin{align}
L_{i+1}&=U_{i+1}^\top M \widetilde{Q}_i \Xi V_{i+1},\label{eq:M_{i+1}}\\
U_{p,i+1}&=(I-U_{i+1}U_{i+1}^\top M) \widetilde{Q}_i\Xi V_{i+1},\label{eq:U_{p,i+1}}\\
V_{p,i+1}&=(I-V_{i+1}V_{i+1}^\top \Xi)\widetilde{Q}_i^\top M U_{i+1}.\label{eq:V_{p,i+1}}
\end{align}
Although \eqref{eq:M_{i+1}}--\eqref{eq:V_{p,i+1}} appear initially unrelated to \eqref{eq:factorM}--\eqref{eq:factorVp}, they are actually closely connected. To see this, we expand $U_{i+1}=U_{i}+\Delta t\dt U_i +\mathcal{O}(\Delta t^2)$ and similarly for $V_{i+1}$. Substituting these expansions into \eqref{eq:M_{i+1}} we get
\begin{equation*}
\begin{aligned}
L_{i+1}&=\left(U_{i}+\Delta t\dt U_i +\mathcal{O}(\Delta t^2)\right)^\top M \left(Q_i+\Delta t G(t_i,Q_i)\right)\Xi \left(V_{i}+\Delta t\dt V_i +\mathcal{O}(\Delta t^2)\right)\\
&=M_i+\Delta t \left(U_i^\top M G(t_i,Q_i)\Xi V_i + U_i M Q_i \Xi \dt V_i + \dt U_i^\top M Q_i \Xi V_i\right)+\mathcal{O}(\Delta t^2) \\
&=M_i+\Delta t \left(U_i^\top M G(t_i,Q_i)\Xi V_i + V_{p,i}^\top \Xi \dt V_{i} + \dt U_i^\top M U_{p,i}\right) +\mathcal{O}(\Delta t^2),
\end{aligned}
\end{equation*}
that is, \eqref{eq:M_{i+1}} approximates \eqref{eq:factorM} up to a second--order error in $\Delta t$. The same argument can be repeated for \eqref{eq:factorUp} and \eqref{eq:factorVp}.
The next lemma formalizes the analysis of the local truncation error and casts light on the dependence on the local curvature.
Throughout the remainder of this subsection, the reference curve $Y$ is regarded as given exactly (i.e., $Y_i=Y(t_i)$ for every $i$). The effect of replacing $Y(t_i)$ by a discrete DLRA approximation is discussed in Remark~\ref{rmk:discrete-Y} below.

\begin{lemma}[Local consistency of \eqref{eq:projected_Euler}]
\label{lem:adj-consistency}
Let $Y\in C^2(\overline{I};\cMr)$, $Q\in C^2(\overline{I};\Rmn)$ be the exact solution of \eqref{eq:covariant_ODE}, and let $Q_{i+1}$ be generated by one step of \eqref{eq:projected_Euler} starting from $Q_i=Q(t_i)$, with $Y_{i+1}=Y(t_{i+1})$, $i\in\{0,\dots,N_t-1\}$. Then, the projected Euler scheme satisfies the local truncation estimate
\begin{equation*}
  \Bigl\| Q(t_{i+1}) - Q_{i+1} \Bigr\|_{M,\Xi}
  \;\le\; C \,\Delta t^2 ,\quad \text{as }\Delta t\rightarrow 0,
\end{equation*}
where 
\begin{equation*}
\begin{aligned}
&C:=\max\limits_{t\in \overline{I}}\kappaeff(t)\|\dt Q(t)\|_{M,\Xi} + \frac{1}{2}\max\limits_{t\in \overline{I}} \|\ddt Q(t)\|_{M,\Xi} +\frac{1}{2}\max\limits_{t\in \overline{I}}\Bigl\|\frac{\rd^2}{\rd t^2}\cP_{Y(t)}\Bigr\|_{\op}\,\max\limits_{t\in \overline{I}}\left(\|Q(t)\|_{M,\Xi}+T \|\dt Q(t)\|_{M,\Xi}\right).
\end{aligned}
\end{equation*}
\end{lemma}
\begin{proof}
We start by considering the Taylor expansions of the maps $\cP: t\mapsto \cP_{Y(t)}$ and $Q: t\mapsto Q(t)$, 
\begin{equation*}
\begin{aligned}
\cP_{Y(t_{i+1})}&=\cP_{Y(t_{i})}+\Delta t D\cP_{Y(t_{i})}\left[\dt Y(t_{i})\right]+ \Delta t^2 R^{\cP},\quad
R^{\cP}:=\int_0^1 (1-s)\,\frac{\rd^2}{\rd t^2}\cP_{Y(t)}\Big|_{t=t_i+s\Delta t}\rd s,\\
Q(t_{i+1})&=Q(t_{i})+\Delta t \dt Q(t_{i})+\Delta t^2 R^Q,\quad R^Q:=\int_0^1 (1-s)\ddt Q(t_i+s\Delta t)\rd s.
\end{aligned}
\end{equation*}
Using that $Y(t_i)=Y_i$ for every $i$, we have that
\begin{equation*}
\begin{aligned}
 Q(t_{i+1}) - Q_{i+1} &= Q(t_{i+1}) - \cP_{Y(t_{i+1})}\left(Q(t_i)+\Delta t \cP_{Y_{i}}\dt Q(t_i)\right) \\
&= Q(t_{i+1})- \left(\cP_{Y(t_i)}+\Delta t D\cP_{Y(t_i)}[\dt Y(t_i)]+\Delta t^2 R^{\cP}\right)\left(Q(t_i)+\Delta t \cP_{Y(t_i)}\dt Q(t_i)\right)\\
&=\Delta t^2\left(R^Q -  D\cP_{Y(t_i)}[\dt Y(t_i)]\cP_{Y(t_i)}\dt Q(t_i) -R^{\cP}Q(t_i)-\Delta t R^{\cP}\cP_{Y(t_i)}\dt Q(t_i) \right),
\end{aligned}
\end{equation*}
where in the last step we used \eqref{eq:derivativetimeQ} with $t=t_i$. Taking the $\|\cdot\|_{M,\Xi}$-norm on both sides and using the triangle inequality leads to the claim.
\end{proof}

\begin{theorem}[Global first-order convergence of the projected Euler method]
\label{thm:global-convergence}
Under the assumptions of Lemma~\ref{lem:adj-consistency}, and assuming moreover that the map $G$ is Lipschitz continuous in its second argument with constant $L_G$ (uniformly in $t$), and that \eqref{eq:projected_Euler} is started from the initial value $Q_0$ of \eqref{eq:covariant_ODE}, the projected Euler method satisfies the global error bound
\begin{equation*}
\|Q(t_i)-Q_i\|_{M,\Xi}\le C\,T\,e^{L_G T}\,\Delta t
\end{equation*}
for all $i\in\left\{0,\dots,N_t\right\}$, where $C$ is the constant of Lemma~\ref{lem:adj-consistency} (hence depends on bounds for $Q$, $\dt Q$, $\ddt Q$, $\kappaeff$, and the second derivative of the projector) and is independent of $\Delta t$.
\end{theorem}

\begin{proof}
Let $e_i:=Q(t_i)-Q_i$ and let $\Phi(t,Q)=\mathcal{P}_{Y(t+\Delta t)}(Q+\Delta t\,G(t,Q))$ denote one step of the method (with exact $Y$). Define the local truncation error $\delta_{i+1}=Q(t_{i+1})-\Phi(t_i,Q(t_i))$ which satisfies $\|\delta_{i+1}\|_{M,\Xi}\le C\Delta t^2$ by Lemma~\ref{lem:adj-consistency}.

Since the orthogonal projector is nonexpansive and $G$ is Lipschitz, the numerical step is Lipschitz:
\[
\|\Phi(t_i,Q)-\Phi(t_i,\widetilde{Q})\|_{M,\Xi}\le(1+L_G\Delta t)\|Q-\widetilde{Q}\|_{M,\Xi}.
\]
Splitting the global error with $e_{i+1} = Q(t_{i+1}) \pm \Phi(t_i, Q(t_i)) - \Phi(t_i, Q_i)$ therefore yields the recurrence
\[
\|e_{i+1}\|_{M,\Xi}\le (1+L_G\Delta t)\|e_i\|_{M,\Xi}+C\Delta t^2.
\]
Unrolling from $e_0=0$ and summing the resulting geometric series gives
\[
\|e_i\|_{M,\Xi}\le C\Delta t^2\sum_{j=0}^{i-1}(1+L_G\Delta t)^{i-1-j}\le C\Delta t\cdot\frac{(1+L_G\Delta t)^i-1}{L_G}\le C \,T\,e^{L_G T}\,\Delta t,
\]
which is the claimed bound. 
\end{proof}

\begin{remark}[Using a discrete state trajectory]\label{rmk:discrete-Y}
Lemma~\ref{lem:adj-consistency} and Theorem~\ref{thm:global-convergence} are proved under the assumption that the projectors are those of the exact curve \(Y(t)\). In computations one only has a discrete trajectory \(\{Y_i\}\) produced by, e.g., a robust DLRA integrator.
The error theory for those integrators \cite{lubichoseledets2014,kieri2019projection} states that, whenever the normal component of the vector field is bounded by \(\varepsilon_Y\),
\[
\|Y_i-Y(t_i)\|_{M,\Xi}\le C(\delta+\varepsilon_Y+\Delta t)
\]
with a constant independent of the singular values. The same theory also yields exactness: if the normal component vanishes, the discrete trajectory coincides with the exact rank-\(r\) solution.

A fully rigorous transfer of this robustness to the difference between one projected-Euler step that uses \(\mathcal{P}_{Y_i}\) and one that uses \(\mathcal{P}_{Y(t_i)}\) would require repeating the structural arguments of the state-error analysis (exactness on the pure rank-\(r\) evolution and avoidance of a direct Lipschitz estimate for the projector). We do not carry out that calculation here. Nevertheless, numerical experiments strongly indicate that the additional consistency error remains of order \(O(\varepsilon_Y+\Delta t)\) and does not introduce a step-size restriction linked to small singular values. A detailed verification is left for future work.
\end{remark}

\paragraph{\textbf{Numerical experiment}}
We consider the time dependent matrix $A(t)=U(t)S(t)V(t)^\top=e^{tW_1}e^t D (e^{tW_2})^\top$, evolving in $\mathbb{R}^{N\times N}$ with $N=100$. $D$ is diagonal with entries $d_j=2^{-j}$ and $W_1$ and $W_2$ are two antisymmetric matrices. The corresponding differential equation is
$$\dt{A}(t) = W_1 A(t) + A(t) + A(t) W_2, \quad A(0) = D.$$
This is an extensively used test case for DLRA integrators since the singular values of $A$ tend to zero exponentially fast \cite{kochlubich2007,kieri2016robust,kieri2019projection,ceruti2022unconventional}. Here, we take
\begin{equation}
W_1=\begin{pmatrix}
0 & -B_1\\
B_1^\top & C_1
\end{pmatrix},\quad W_2=\begin{pmatrix}
0 & -B_2\\
B_2^\top & C_2
\end{pmatrix},
\end{equation}
where $B_j$ are random matrices of dimensions $n_1\times (N-n_1)$ and $C_j$ are random antisymmetric matrices of dimensions $(N-n_1)\times (N-n_1)$ for $j=1,2$. Denoting with $E_r\in \mathbb{R}^{N\times r}$ the matrix containing the first $r$ canonical basis vectors, the best rank-$r$ approximation of $A(t)$ at each time instant is given by $A_r(t)=U_r(t)S_r(t)V_r(t)^\top$, where $U_r(t)=U(t)E_r$, $S_r(t)=E_r^\top S(t)E_r$ and $V_r(t)=V(t)E_r$. The chosen structure of $W_1$ and $W_2$ guarantees that $\dt U_r(t)^\top U_r(t)=\dt V_r(t)^\top V_r(t)=0$ for any $r\in\{1,\dots,n_1\}$.
We then consider a time-dependent reference vector field $Q(t)$ characterized by $L(t)=\exp(t)I_{r}$, $U_p(t)=U(t)E^\perp_r$ and $V_p(t)=V(t)E^\perp_r$, $E^\perp_r$ collecting the canonical vectors ${\bm e}_{r+1},\dots, {\bm e}_{2r}$.
The left panel of Figure \ref{fig:NumTestCovariant} shows the error at the final time as $\Delta t$ decreases and for $r\in\{4,8,16,32\}$ between the reference vector field and the two approximations obtained either by integrating the factor equations \eqref{eq:factorM}--\eqref{eq:factorVp} with an explicit Runge--Kutta method of order 4, or by the projected Euler scheme \eqref{eq:M_{i+1}}--\eqref{eq:U_{p,i+1}}. We remark that neither approach requires a restriction on the time step despite the increasingly smaller singular values of $A_r(t)$. As shown in Theorem \ref{thm:global-convergence}, the error depends on the effective curvature, which remains of moderate size as $r$ increases (right panel of Figure \ref{fig:NumTestCovariant}).
Note that in this experiment we set $M=\Xi=I_{N}$, $I_N$ being the identity matrix of $\mathbb{R}^{N\times N}$ and $n_1=50$.

\subsubsection{Application to the adjoint equation}\label{subsec:adjoint}
Coming back to our optimization context, \eqref{eq:opt-DLRAPadjoint} is equivalent to a backward formulation of \eqref{eq:covariant_ODE} by setting 
\[G(t,Q(t))=\cP_{Y(t)}\left(F_X(t,Y(t),\zb(t))^\star Q(t)+ X_d(t)-Y(t)\right)+D\cP^\star_{Y(t),F^\perp(t,Y(t),\zb(t))} Q(t).\]
The projected Euler scheme then computes backwards in time
\[Q_{i}=\cP_{Y_i}\left(Q_{i+1}+\Delta t G(t_{i+1},Q_{i+1})\right),\quad i=0,\dots,N_t-1,\]
starting from $Q_{N_t}=0\in T_{Y_{N_t}}\cMr$. 
The computation of $\widetilde{Q}_{i+1}:=Q_{i+1}+\Delta t G(t_{i+1},Q_{i+1})$ can be done efficiently in tangent format by first computing the tangent parametrization $(\widehat{L}_{i+1},\widehat{U}_{p,i+1},\widehat{V}_{p,i+1})$ of $\cP_{Y_{i+1}}\bigl(F_X(t_{i+1},Y_{i+1},\zb_{i+1})^\star Q_{i+1}+ X_d(t_{i+1})-Y_{i+1}\bigr)$. Then, the triplet $(\widetilde{L}_{i+1},\widetilde{U}_{p,i+1},\widetilde{V}_{p,i+1})$ identifying $\widetilde{Q}_{i+1}$ is given by 
\begin{align}
\widetilde{L}_{i+1}&=L_{i+1}+\Delta t\widehat{L}_{i+1},\quad \widetilde{U}_{p,i+1}=U_{p,i+1}+\Delta t\bigl(\widehat{U}_{p,i+1}+ F^\perp_{i+1}\Xi V_{p,i+1} S_{i+1}^{-\top}\bigr),\\
\widetilde{V}_{p,i+1}&=V_{p,i+1} +\Delta t\bigl( \widehat{V}_{p,i+1}+ (F^\perp_{i+1})^\top M U_{p,i+1} S_{i+1}^{-1}\bigr),
\end{align}
where $F^\perp_{i+1}:=(I-\cP_{Y_{i+1}})F(t_{i+1},Y_{i+1},\zb_{i+1})$ and we used \eqref{eq:expression_adjointderivative} to express the curvature term. Note that the curvature term involves the inverse of the singular value matrix of $Y$. This may induce stiffness in the time integration, potentially imposing restrictions on the choice of the step size. We stress that this possible stiffness arises from the particular structure of the right-hand side $G$ in \eqref{eq:covariant_ODE} for our problem under consideration, rather than being a feature of the time integration scheme \eqref{eq:projected_Euler}. In our experiments of Section \ref{sec:num}, we will observe that $F^\perp \Xi V_{p}$ and $U_pMF^\perp$ are typically very small, so that, in practice, no significant restriction on the choice of the time step is observed.

\subsection{Computational cost}\label{sec:computationalcost}
We conclude this section by discussing the computational advantages of the DLRA approaches presented in this work compared to a classical full-rank implementation.
From a memory perspective, at each time instant $t_i$ the system state $Y_i\in \Rmn$ is represented by the triplet $(U_i,S_i,V_i)$, requiring $\mathcal{O}((m+n+r)r)$ storage. Similarly, the adjoint variable $Q_i$ is represented either by the ``SVD'' triplet $(U^Q_i,S^Q_i,V^Q_i)$ in the OtP approach, or by the ``tangent'' triplet $(L_i,U_{p,i},V_{p,i})$ at the bases $(U_i,V_i)$ in the PtO formulation. In both cases, storing the adjoint variable also scales as $\mathcal{O}((m+n+r)r)$. Hence, the total memory footprint of each gradient descent iteration is of the order $\mathcal{O}(N_t(m+n+r)r)$, which can be substantially smaller than the
$\mathcal{O}(N_t mn)$ requirement of a classical full-rank forward-backward iteration.

On the other hand, computational efficiency requires the ability to evaluate both $F$ and its linearized adjoint $F_X ^\star$ directly in a low-rank format, without ever assembling any full dense $(m\times n)$ matrix. Prototype examples are vector fields of the form
\begin{equation}\label{eq:affine_expansion}
F(t,X,\zb)= \mathcal{A}(X)+\fb(t) + B\zb(t),\quad \mathcal{A}(X):=\sum_{k=1}^{N_p} A_k X D_k,
\end{equation}
where $A_k\in\mathbb{R}^{m\times m}$, $D_k\in\mathbb{R}^{n\times n}$, $B\in\mathcal{L}(\Rs,\mathbb{R}^{m\times n})$, $\zb(t)\in \Rs$ and $\fb(t)\in \Rmn$ at each instant of time. In such cases, $\mathcal{A}(X)$ can be assembled in $\mathcal{O}(N_p(m+n)r^2)$
for a given $X\in\cMr$ (assuming $A_k$ and $D_k$ are sparse). Provided that also the action of the control operator can be computed at a comparable cost, the DLRA approaches yield substantial computational savings compared with integrating the full-rank ambient matrix ODE which would require $\mathcal{O}(N_pmn)$ operations. The affine expansion \eqref{eq:affine_expansion} arises in many applications, including the optimal control of parabolic PDEs posed on Cartesian meshes, or the optimal control of random PDEs due to a truncated Karhunen--Lo\`eve representation of the random field~\cite{lord2014introduction} (see Section \ref{sec:ocp}). Nevertheless, it is not strictly required. Nonaffine right-hand sides can be accommodated by interpolatory DLRA (via DEIM and variants) as recently proposed in~\cite{carrel2024interpolatory}.
Table~\ref{tab:cost-model} summarizes the total number of operations and memory footprint of one gradient step in the ambient space, and with the two DLRA approaches.

\begin{table}[t]
  \centering
  \caption{Per gradient iteration cost model: $N_t$ is the number of time steps,
           $m$, $n$ are the matrix dimensions,
           $r$ is the rank, $N_p$ is the number of terms in the expansion \eqref{eq:affine_expansion}.}
  \label{tab:cost-model}
  \begin{tabular}{lll}
    \toprule
    Operation & Ambient Space & DLRA Approaches \\
    \midrule
    Forward sweep
      & $\mathcal{O}\bigl( N_p m n\, N_t \bigr)$
      & $\mathcal{O}\bigl( N_p(m{+}n)\, r^2 N_t \bigr)$ \\
    Backward sweep
      & $\mathcal{O}\bigl( N_p\, m n\, N_t \bigr)$
      & $\mathcal{O}\bigl( N_p(m{+}n)\, r^2 N_t \bigr)$ \\
    \midrule
    Memory footprint
      & $\mathcal{O}\bigl( m n\, N_t \bigr)$
      & $\mathcal{O}\bigl( (m+n+r)\, r\, N_t \bigr)$ \\
    \bottomrule
  \end{tabular}
\end{table}

\section{Application to parabolic optimal control problems with random coefficients}
\label{sec:ocp}

Let $D \subset \mathbb{R}^2$ be a bounded Lipschitz domain, $I\bydefeq(0,T)$ a finite time interval, $(\Omega, \mathcal{F}, \mathbb{P})$ a complete probability space, and $V\bydefeq H^1_0(D)$ the standard Sobolev space.
We consider a physical system described by the function-valued random variable $y:\Omega\rightarrow L^2(I;V)\cap H^1(I;V^\star)=:\mathcal{V}$, solution of the random parabolic equation
\begin{equation}
  \label{eq:pde-state}
  \partial_t y(x, t; \omega)
  \;-\; \nabla_x \!\cdot\! \bigl( a(x, \omega)\, \nabla_x y(x, t; \omega) \bigr)
  \;=\; f(x, t) \;+\; z(x, t)
  \qquad \text{for a.e. (x,t) in } D \times I,\; \PP\text{-a.e.}\;\omega,
\end{equation}
with homogeneous Dirichlet boundary conditions on $\partial D$, and
$y(\cdot, 0; \omega) = y_0(\cdot,\omega)$ as initial condition for $\PP\text{-a.e.}\;\omega$. Here, $a\in L^2(\Omega;L^\infty(D))$ is a random
field, uniformly bounded from above and below, and $f\in L^2(I;L^2(D))$ is a deterministic source. Our goal is to find a deterministic control $z\in \mathcal{Z}\bydefeq L^2\!\bigl(I;\, L^2(D) \bigr)$, that is, a single open-loop control independent of the random realizations. This is the case of interest in applications where design choices have to be made before the randomness is observed. To determine $z$, we consider the cost functional
\begin{equation}
  \label{eq:cost-pde}
  \mathcal{J}(y, z)
  \bydefeq \tfrac{1}{2}\, \mathbb{E}\!\left[
          \int_0^T \!\! \| y(\cdot, t;\, \omega) - y_d(\cdot, t) \|_{L^2(D)}^2 \, \rd t
        \right]
      \;+\; \tfrac{\nu}{2} \int_0^T \!\! \| z(\cdot, t) \|_{L^2(D)}^2 \, \rd t,
\end{equation}
where $y_d \in L^2(I; L^2(D))$ is a target trajectory. The sought control is computed as the solution of
\begin{equation}
  \label{eq:ocp-pde}
  \min_{y \in \mathcal{V}, z \in \mathcal{Z}} \mathcal{J}(y, z)
  \quad \text{subject to \eqref{eq:pde-state}.}
\end{equation}

Under the above assumptions, the state equation admits a unique solution $y(z)\in \mathcal{V}$ for every $z\in \mathcal{Z}$, and the reduced functional $\widehat{\mathcal{J}}: \mathcal{Z}\rightarrow \mathbb{R}$ is defined as $\widehat{\mathcal{J}} (z) = \cJ(y(z),z)$. Since the map $z\mapsto y(z)$ is affine, it is immediate to verify that $\widehat{\mathcal{J}}$ is convex, and the existence of a unique minimizer $z^\star$ of \eqref{eq:ocp-pde} follows from classical arguments of calculus of variations, see, e.g., \cite[Theorem 1]{kouri2018optimization} and \cite[Section 3]{guth2024parabolic} for a related discussion.

\subsection{Semidiscretization and optimality conditions}
\label{sec:mc-semidisc}

To numerically approximate the solution of \eqref{eq:ocp-pde},
we replace the expectation operator by a sample-based quadrature formula $\{(\omega_j, w_j)\}_{j=1}^n$, so that for every $g\in L^2(\Omega)$, 
\[\mathbb{E}[g(\omega)] \approx \sum_{j=1}^n w_j\, g(\omega_j).\]
We only consider quadrature formulae with positive
weights, i.e., $w_j > 0$ for $j\in\{1,\dots,n\}$, since negative weights may compromise the convexity of the discretized functional \cite{nobile2024combination,nobile2025multilevel}. This class includes, for instance,
Monte Carlo, quasi-Monte Carlo and Gaussian quadrature methods. The quadrature weights are collected into a diagonal and positive definite matrix
\begin{equation}\label{eq:Xi-def}
\Xi\bydefeq\mathrm{diag}\bigl(w_1,\dots,w_n\bigr),\quad \Xi\succ 0.
\end{equation}
Concerning the spatial discretization, we employ the Finite Element Method (FEM) based on a continuous finite-element space $V_h \subset H^1_0(D)$
of dimension $m$. We denote with $A(\omega_j) \in \R^{m \times m}$, $j\in \{1,\dots,n\}$, the stiffness matrices associated to the realizations $a(\cdot,\omega_j)$ of the random coefficient and with $M$ the mass matrix. The time-dependent vectors $\yb_j(t)\in \mathbb{R}^m$, $j\in\{1,\dots, n\}$ and for $t\in \overline{I}$, contain the degrees of freedom of the $j$-th realization of the state variable $y(\cdot,t;\omega_j)$ expressed in the FE basis. Similarly, $\mathbf{f}(t)\in \mathbb{R}^m$, $\yb_d(t)\in \mathbb{R}^m$, $\zb(t)\in \mathbb{R}^m$ collect the degrees of freedom of the force term, target state and control at each time instant, and $\yb_{0,j}\in\mathbb{R}^m$ of the random initial conditions for $j\in\{1,\dots,n\}$.
Note that, for this example, $s=m$.

With these two semidiscretizations at hand, the semidiscrete approximation of \eqref{eq:ocp-pde} is
\begin{equation}
\begin{aligned}
  \label{eq:ocp-pde_discretized}
  & \min_{\substack{\{\yb_j\}_{j=1}^n \subset H^1(I; \mathbb{R}^m)\\
        \zb\in  L^2(I; \mathbb{R}^m)
        }} \mathcal{J}_{h,n}(\{\yb_j\}_{j=1}^n,\zb) \bydefeq \tfrac{1}{2}\, \sum_{j=1}^n w_j \left[
          \int_0^T \!\! \| \yb_j(t) - \yb_d(t) \|_{L^2(D)}^2 \, \rd t
        \right]
      \;+\; \tfrac{\nu}{2} \int_0^T \!\! \| \zb(t) \|_{L^2(D)}^2 \, \rd t,\\
    & \text{subject to }\;  \begin{cases} 
    &M \dt{\yb}_j(t)   \;+\; A(\omega_j) \yb_j(t) \;=\; M\left(\mathbf{f}(t) \;+\; \zb(t)\right),\quad \text{for a.e. } t\in (0,T],\;j=1,\dots,n,\\
    & \yb_j(0)=\yb_{0,j},\quad\ j=1,\dots,n,
    \end{cases}
 \end{aligned}
\end{equation}
where we used the slight abuse of notation $\|\vb\|^2_{L^2(D)}$ to denote the $L^2(D)$ norm of the finite element function represented by the vector $\vb$ in the FEM basis. The existence and uniqueness of the minimizer of \eqref{eq:ocp-pde_discretized} follows as for \eqref{eq:ocp-pde}.
The optimality conditions are obtained by either setting the Fr\'echet derivative of $\widehat{\mathcal{J}}_{h,n}$ equal to zero, see, e.g., \cite[Section 3]{guth2024parabolic}, or by
setting the partial derivatives of the associated Lagrangian function.
The optimality system reads
\begin{align}
&M \dt{\yb}_j(t)
  \;+\; A(\omega_j) \yb_j(t) \;=\; M\left(\mathbf{f}(t) \;+\; \zb(t)\right),\quad &\text{for a.e. } t\in (0,T],\;j\in\{1,\dots,n\},\label{eq:state_eq_discretized}\\
&-M \dt{\mathbf{p}}_j(t)
  \;+\; A(\omega_j) \mathbf{p}_j(t) \;=\; M\left(\yb_d(t)\;-\;\yb_j(t)\right),\quad &\text{for a.e. } t\in [0,T),\;j\in\{1,\dots,n\},\label{eq:adjoint_eq_discretized}\\
&\nu \zb(t)-\sum_{j=1}^n w_j \mathbf{p}_j(t)=0,\quad &\forall t\in (0,T],\label{eq:optcond_eq_discretized}    
\end{align}
with the additional initial conditions $\yb_j(0)=\yb_{0,j}$ and terminal conditions $\mathbf{p}_j(T)=0$, for $j\in\{1,\dots,n\}$.
Note that in \eqref{eq:adjoint_eq_discretized} we used that the matrices $A(\omega_j)$ are all symmetric. 

\subsection{Affine random expansion and matrix reformulation}
\label{sec:affine}

We assume an affine expansion of the random diffusion coefficient, namely
\begin{equation}
  \label{eq:affine-coefficient}
  a(x,\omega)=a_0(x)+\sum_{k=1}^{N_p} \zeta_k(\omega)a_k(x),
\end{equation}
where $\left\{a_k\right\}_{k=0}^{N_p}$ are deterministic coefficients and $\{\zeta_k\}_{k=1}^{N_p}$ are scalar-valued random variables.
After the discretization, this hypothesis implies that the stiffness matrices can be expressed as
\begin{equation}
  \label{eq:affine-stiffness}
  A(\omega_j) \;=\; A_0 + \sum_{k=1}^{N_p} \zeta_k(\omega_j)\, A_k,\quad j\in\{1,\dots,n\}
\end{equation}
where $A_{k} \in \mathbb{R}^{m \times m}$ is the deterministic
stiffness matrix related to the coefficient $a_k$, $k\in\{0,\dots,N_p\}$.

Next, we represent the ensemble of discrete random states by the time-dependent matrix
\begin{equation}
  \label{eq:Y-matrix}
  X(t)=\left[\yb_1(t)|\dots |\yb_n(t)\right] \in \mathbb{R}^{m \times n},
\end{equation}
whose $j$-th column coincides with the degrees of freedom of the $j$-th state realization.
Exploiting the affine structure \eqref{eq:affine-coefficient}, the state constraint of \eqref{eq:ocp-pde_discretized}
can be formulated as the matrix ODE,
\begin{equation}
  \label{eq:Y-ode}
  \begin{cases}
  \dt X(t) &=\; F\bigl(t, X(t),\, \zb(t) \bigr),\quad \text{for a.e. }t\in (0,T],\\
   X(0) &= X_0,
  \end{cases}
\end{equation}
where the $j$-th column of $X_0\in \mathbb{R}^{m\times n}$ corresponds to $\yb_{0,j}$, and
the right-hand side $F:\mathbb{R}\times \mathbb{R}^{m\times n}\times \mathbb{R}^m\rightarrow \mathbb{R}^{m\times n}$ is defined as
\begin{equation}
  \label{eq:F-def}
  F(t, X, \zb)
  \bydefeq -\, \mathcal{A}(X) \;+\; B\,(\zb + \mathbf{f}(t)),
  \qquad \mathcal{A}(X) \bydefeq M^{-1} \!\!\sum_{k=0}^{N_p}\!\! A_k\, X\, D_k.
\end{equation}
Here, $D_0=I_n$ is the identity matrix of $\mathbb{R}^n$ and, for $k\in \{1,\dots,N_p\}$, $D_k$ are diagonal matrices collecting the sample-dependent realizations of the $k$-th random variable,
\begin{equation}\label{eq:definition_Dk}
D_k=\mathrm{diag}\bigl(\zeta_k(\omega_1),\dots,\zeta_k(\omega_n)\bigr),
\end{equation}
and
$B :\mathbb{R}^m\rightarrow \mathbb{R}^{m\times n}$, $B\vb \bydefeq \vb\, \mathbf{1}_n^\top$ is the
column-broadcast operator. Note that $\mathcal{A}$ is self-adjoint with respect to $\langle \cdot, \cdot\rangle_{M,\Xi}$, since each $A_k$ is symmetric and each $D_k$ is diagonal.
Finally, due to the definition of the $\|\cdot\|_{M,\Xi}$ norm, we have that
\[\frac{1}{2}\int_0^T \|X(t)-X_d(t)\|^2_{M,\Xi}\rd t=\frac{1}{2}\sum_{j=1}^n w_j \left[\int_0^T \|\yb_j(t)-\yb_d(t)\|_{L^2(D)}\rd t\right],\]
with $X_d(t):=\yb_d(t)\mathbf{1}_n^\top\in\Rmn$, so that \eqref{eq:ocp-pde_discretized} can be recast in the form of \eqref{eq:ocp_full}.

\subsection{Numerical experiments}
\label{sec:num}
In this subsection we present numerical results aimed at testing the OtP and PtO approaches applied to the minimization of \eqref{eq:ocp-pde_discretized}, and at showing their computational advantages compared to a classical full-rank approach. Additionally, we investigate a few geometric properties such as the time evolution of $\kappaeff$ and of the bounds appearing in Assumptions \ref{ass:modelling_error_F}, \ref{ass:modelling_error_Fx} and \ref{ass:modelling_error_Fx2}. The experiments have been performed in \textsc{Matlab} 2023b. The codes and data analyzed in this section are available at \cite{GitHub}.

We set $\mathcal{D}=(0,1)^2$, $T=1$ and $f=0$. The spatial discretization is performed using $\mathbb{P}^1$ finite elements on a uniform mesh. We consider the random diffusion coefficient \eqref{eq:affine-coefficient}, with
\[a_0(x)=\bar{a},\quad a_k(x)=\frac{\bar{a}}{(j_k^2+\ell_k^2)}\sin(j_k\pi x_1)\sin(\ell_k\pi x_2),\]
where the sequence $(j_k,\ell_k)_{k\in\mathbb{N}}$ is an ordering of the elements of $\mathbb{N}\times\mathbb{N}$ and $\bar{a}:=10^{-2}$. We truncate the expansion at $k=4$, including only the pairs $(1,1)$, $(1,2)$, $(2,2)$ and $(2,3)$. The random variables $\{\zeta_k\}$ are uniformly distributed over $(0,1)$. The discretization of the probability space is performed with the Monte Carlo method. The matrix $\Xi$ is thus diagonal with equal weights $1/n$, $n$ being the number of samples.

The interval $(0,T]$ is divided into $N_t\in \N$ subintervals $I_i:=(t_i,t_{i+1}]$, $i=0,\dots,N_t-1$, of equal length $\Delta t$.
To discretize in time the full-rank optimization problem, we follow a discretize-then-optimize paradigm: the objective functional is discretized by approximating the time integrals with a left-endpoint rectangular quadrature rule, while the state equation \eqref{eq:Y-ode} is discretized using the explicit Euler method. The adjoint equation and the optimality condition are then discretized consistently so as to obtain an \textit{exact} discrete gradient.
The time discretization for both the OtP and PtO approaches is detailed in Section \ref{sec:disc}.
To speed up computations, we have used consistently the lumped mass matrix \cite[\S 11.4]{quarteroni1994numerical} for both the time integration and as matrix $M$ in the weighted Frobenius norm.
For all experiments, the target state is the deterministic function $y_d(x,t)=t\sin(2\pi x_1)\sin(2\pi x_2)$ and the regularization parameter is set to $\nu=10^{-3}$. As random initial state we consider a Gaussian process characterized by a zero mean and by the covariance function $\text{Cov}(x,x'):=\sigma^2\exp(-\|x-x'\|_2^2)$ with $\sigma^2=0.1$. The resulting matrix $X_0$ has fast decaying singular values.

\begin{figure}
\centering
\includegraphics[scale=0.4]{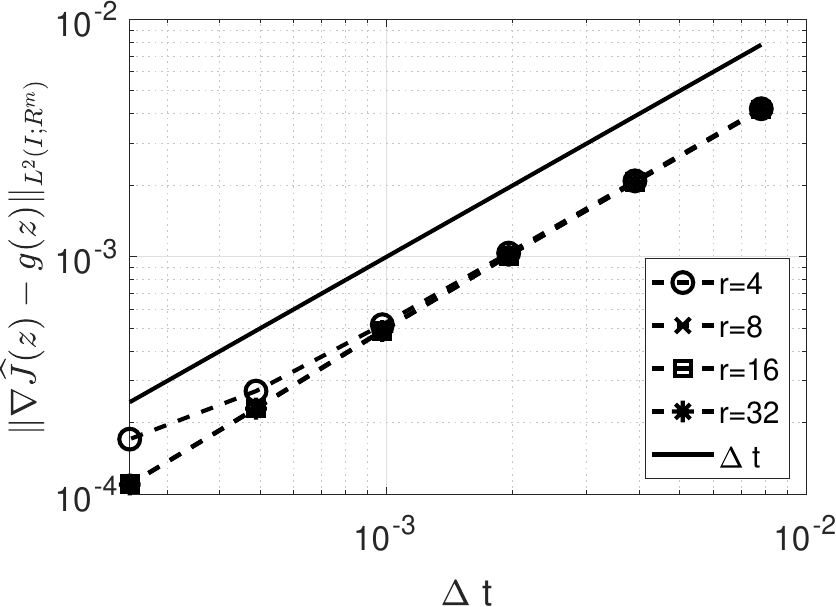}
\includegraphics[scale=0.4]{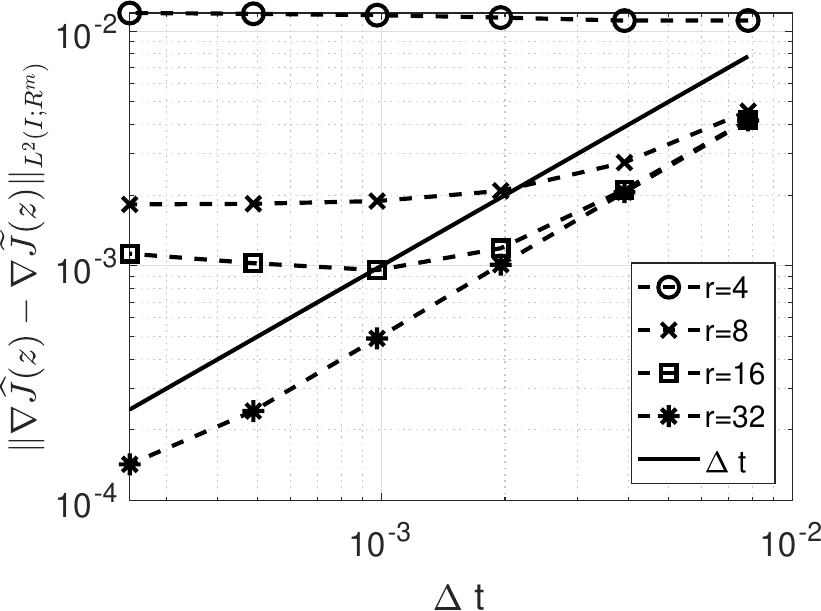}
\includegraphics[scale=0.4]{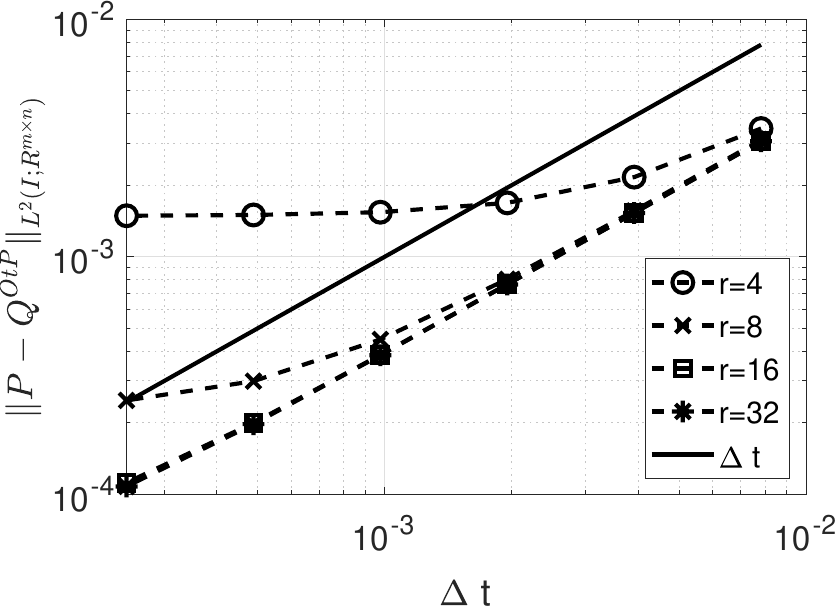}
\includegraphics[scale=0.4]{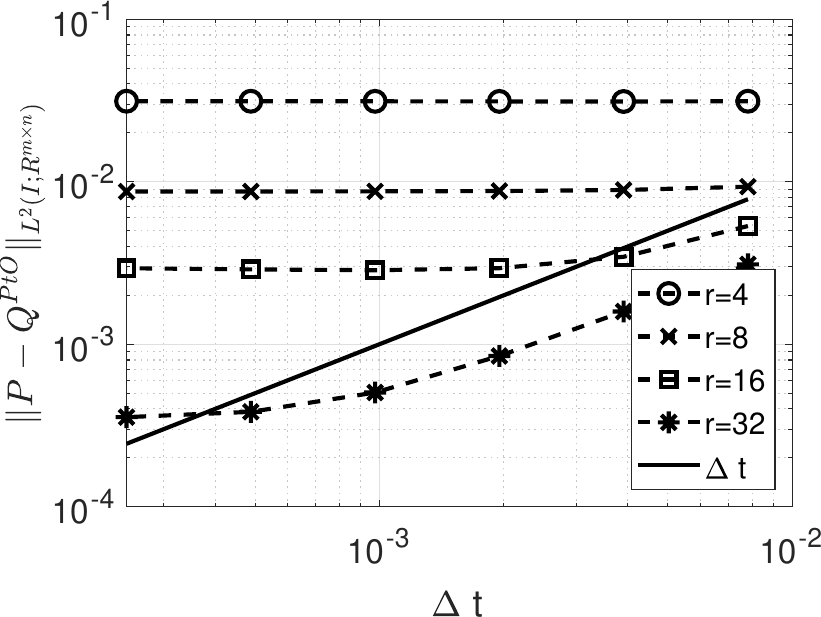}
\caption{For different values of $\Delta t$, we report the error in the discretized $L^2(I;L^2(\mathcal{D}))$ norm between the full-rank gradient (top row) and full-rank adjoints (bottom row) with their OtP (left column) and PtO (right column) approximations. Parameters: $m=961$, $n=500$.}
\label{fig:gradients}
\end{figure}

Figure \ref{fig:gradients} shows the error, measured in the discretized $L^2(I;L^2(D))$-norm and for different values of $\Delta t$, between the full-rank adjoint and gradient, and their corresponding OtP and PtO approximations. The gradients and adjoints are computed starting from $z(x,t)=t\exp(x_1+x_2)\sin(2\pi x_2)$. Both DLRA approaches converge with a rate $\Delta t$ before the truncation error dominates. The PtO paradigm exhibits, however, a larger error for given $\Delta t$ and $r$. This is not surprising since the OtP approach computes a ``true'' dynamical low-rank approximation of the adjoint variable, while the PtO computes a low-rank approximation by projecting the dynamics onto the tangent space at the given state variable. 
Concerning the time integration, the PtO approach is not subject to any particular restriction on the time-step choice, as anticipated in Section \ref{subsec:adjoint}. For these parameter values, the CFL condition for the full-rank ODE is approximately $1.6\cdot 10^{-2}$.
Finally, we emphasize that the gradient approximations benefit from the action of the adjoint control operator $B^\star$, which maps the matrix of adjoint states onto an aggregated quantity, that is, their mean. This averaging procedure mitigates approximation errors, leading to gradient approximations that are typically more accurate than the corresponding adjoint approximations.

\begin{figure}
\centering
\includegraphics[scale=0.4]{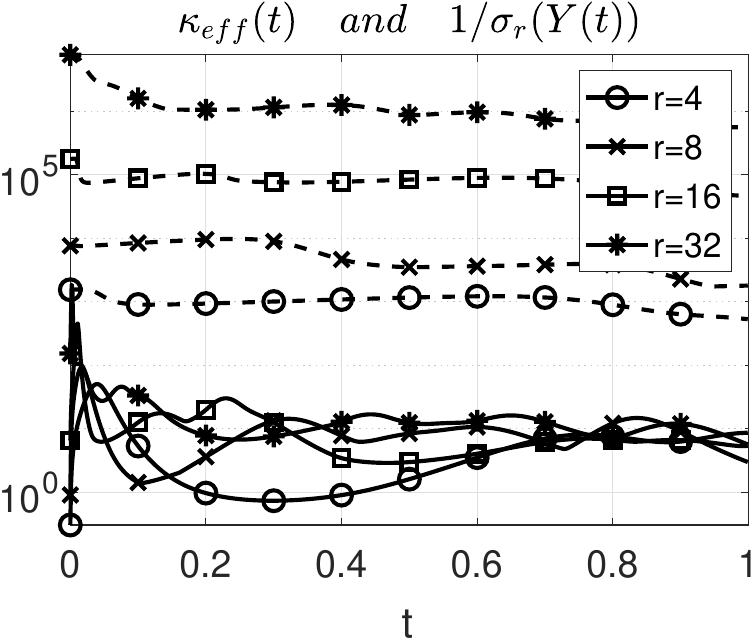}
\includegraphics[scale=0.4]{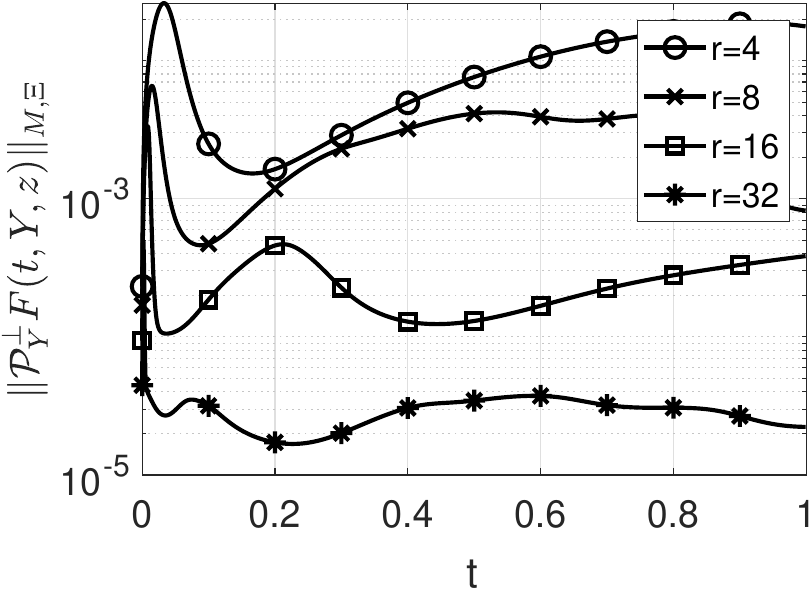}\\
\includegraphics[scale=0.4]{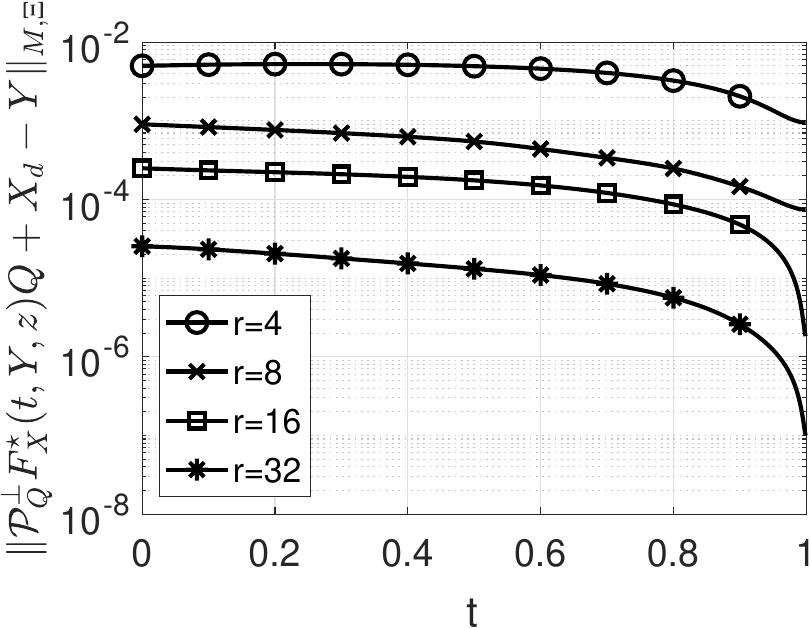}
\includegraphics[scale=0.4]{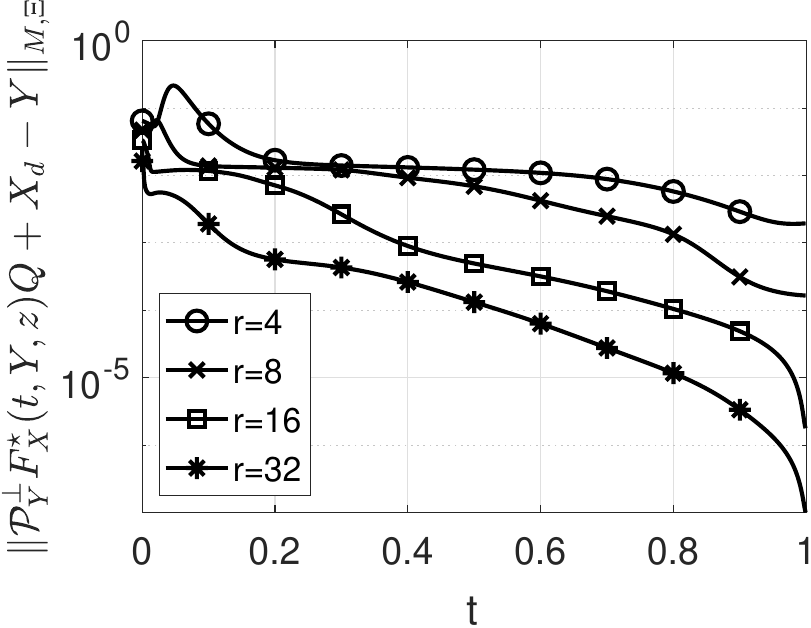}
\caption{Time evolution of the effective curvature (continuous line) and of $1/\sigma_r(Y(t))$ (dashed line) in the top-left panel, and of the normal components of the state vector field (top-right), of the adjoint vector field for the OtP approach (bottom-left) and for the PtO approach (bottom-right). Parameters: $m=961$, $n=500$, $\Delta t=10^{-3}$ and control $z(x,t)=t\exp(x_1+x_2)\sin(2\pi x_2)$.
}\label{fig:geometricanalysis}
\end{figure}

The top left panel of Figure \ref{fig:geometricanalysis} shows the time evolution of the effective curvature $\kappaeff(t):=D\cP_{Y(t)}[\dt Y(t)]$ and of $1/\sigma_r(Y(t))$. The effective curvature is estimated through the bound
\begin{equation*}
\begin{aligned}
&\kappaeff:=\|D\cP_{Y}[\dt Y]\|_{\op}\leq \sup_{\|Z\|_{M,\Xi}=1} \left(\|DP_U[\dt Y]Z(P_V^\perp)^\top\|_{M,\Xi}+\|P_U^\perp Z \left(DP_V[\dt Y]\right)^\top\|_{M,\Xi}\right)\\
&\leq \sup_{\|Z\|_{M,\Xi}=1}\left(\|M^{\frac{1}{2}} DP_U[\dt Y]M^{-\frac{1}{2}}\|_2\|Z(P_V^\perp)^\top\|_{M,\Xi}
+ \|P_U^\perp Z\|_{M,\Xi}\|\Xi^{\frac{1}{2}} DP_V[\dt Y]\Xi^{-\frac{1}{2}}\|_2\right)\\
&\leq \|M^{\frac{1}{2}} DP_U[\dt Y]M^{-\frac{1}{2}}\|_2 +\|\Xi^{\frac{1}{2}}DP_V[\dt Y]\Xi^{-\frac{1}{2}}\|_2,
\end{aligned}
\end{equation*}
which is evaluated using \eqref{eq:directional_derivatives_PU_PV} and \eqref{eq:dotU_dotV}. We emphasize that the effective curvature does not explode as $\sigma_r(Y(t))\rightarrow 0$, and remains well-behaved as the rank increases.
The top right panel of Figure \ref{fig:geometricanalysis} shows the time evolution of the normal component of the vector field $F$ for different rank values, and similarly, the bottom panels show the evolution of the normal component of the adjoint vector fields projected either on $Q$ (left panel) or on the current state $Y$ (right panel). Recall that the adjoint equations are integrated backwards in time, and this explains why the normal components are growing as $t$ decreases. It is evident that the PtO approach has a considerably larger normal component, and this explains the larger error constants observed in Figure \ref{fig:gradients}.

\begin{figure}
\centering
\includegraphics[scale=0.45]{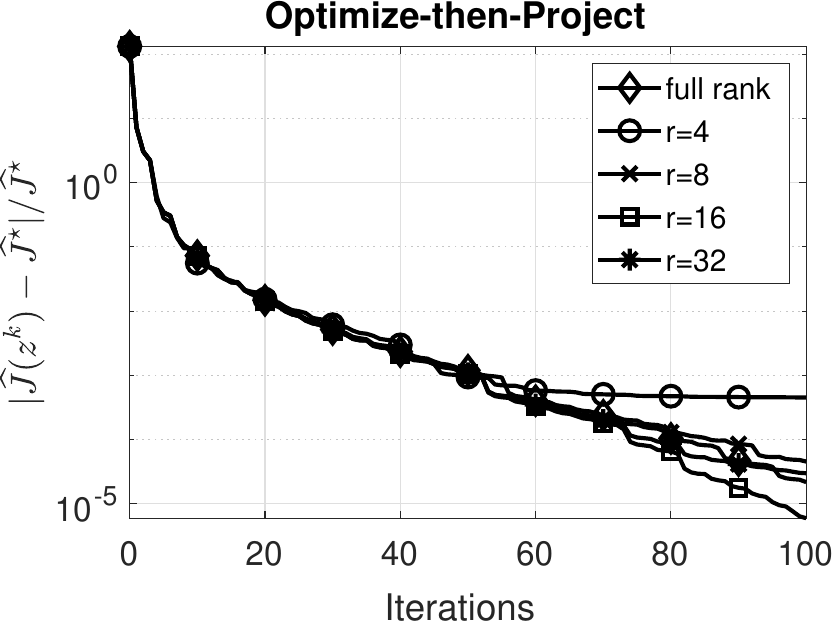}
\includegraphics[scale=0.45]{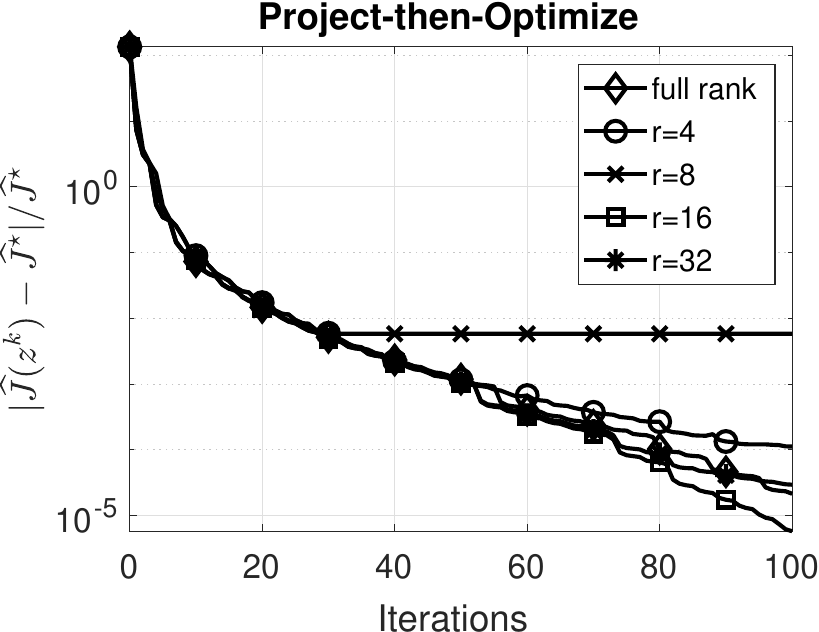}
\caption{Evolution of the objective function value gap during gradient descent iterations for different rank values. Markers are shown every 10 iterations to improve visualization. Parameters: $m=961$, $n=500$, $T=1$, $\Delta t=10^{-3}$.}
\label{fig:optimization}
\end{figure}

Figure \ref{fig:optimization} shows the relative objective function value gap, defined as the absolute difference between a reference value and the value of the full-rank objective functional evaluated on the iterates of the gradient descent on the approximated problems divided by the reference value, for different values of the rank $r$. The reference value $\widehat{\cJ}^\star$ is computed by 
running the gradient descent method on the full-rank problem up to a tolerance of $10^{-6}$ on the $L^2(I;L^2(D))$-norm of the gradient. In total, $163$ iterations are needed to reach this accuracy.
We observe that even small values of the rank are sufficient to achieve quite small objective function value gaps. Note that, for the PtO approach with $r=8$, the optimization algorithm appears to get stuck in a local minimum (gradient descent requires about $8$ line search steps, while for all other rank values the average number is between $2$ and $3$), confirming that the PtO approach may lead to complex optimization landscapes, even though the original problem is convex (linear-quadratic for this example), see Remark \ref{rmk:nonconvexity}. On the left panel, we observe instead that with $r=16$ we achieve even better accuracies than the full-rank approach after the same number of iterations.

\begin{figure}[htbp]
\centering
\begin{minipage}[c]{0.47\textwidth}
    \centering
    \includegraphics[width=0.8\linewidth]{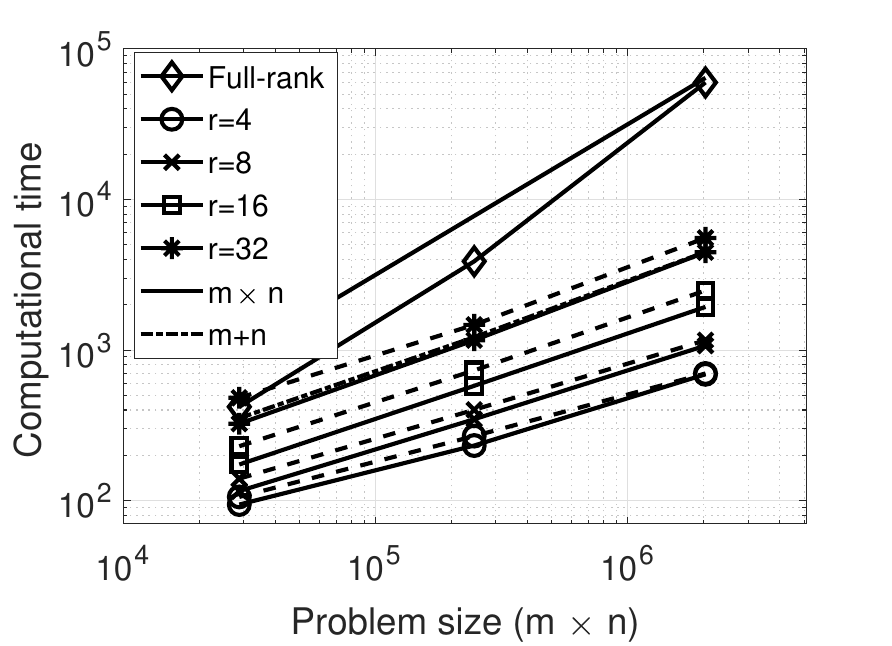}
\end{minipage}\quad
\begin{minipage}[c]{0.50\textwidth}
    \centering
    \begin{tabular}{cccccc}
    \toprule
    $(m, n)$ & $r = 4$ & $r = 8$ & $r = 16$ & $r = 32$ & Full\\
    \midrule
    $225 \times 128$   & 0.02 & 0.05 & 0.09 & 0.20 & 0.46\\
    $961 \times 256$   & 0.08 & 0.16 & 0.32 & 0.64 & 3.94\\
    $3969 \times 512$  & 0.29 & 0.58 & 1.15 & 2.31 & 32.55\\
    \bottomrule
    \end{tabular}
\end{minipage}
\caption{On the left, computational times in seconds of 30 gradient descent iterations for the full-rank, OtP (continuous line) and PtO (dashed line) approaches for different rank values. On the right, comparison of the memory footprints expressed in gigabytes to store the state and adjoint variables. Parameters: $m\in \{225,961,3969\}$, $n\in\{128,256,512\}$, $T=1$, $\Delta t=10^{-3}$.}
\label{fig:TimeOptimization}
\end{figure}

Finally, Figure~\ref{fig:TimeOptimization} reports the computational times required by the full-rank, the OtP and PtO approaches to perform 30 iterations of the gradient descent method. We consider three problem sizes obtained by discretizing in space on a uniform grid of mesh sizes $2^{-\ell}$, $\ell\in\{4,5,6\}$, resulting in $m\in \{225,961,3969\}$, and with $n\in\{128,256,512\}$ samples.
The complexity of the full-rank approach scales as $\mathcal{O}(m n)$ as predicted, while that of both OtP and PtO is of the order $\mathcal{O}((m+n)r)$. The table on the right compares the memory required to store the state and adjoint variables for the dynamical low-rank and the full-rank approaches. Both OtP and PtO require the same memory since their adjoints are stored as three factors of identical dimensions.
The speed-up and memory savings are impressive and highlight the efficiency of the computational strategies presented in this work for this class of optimal control problems.

\section{Conclusion}
\label{sec:conclusion}
We have presented a mathematical framework to use dynamical low-rank approximation strategies to accelerate the optimal control of systems described by time-dependent matrices. Our analysis casts light on the interplay between the magnitude of normal components of vector fields, system stability, and geometric factors, such as the effective curvature, which collectively determine the accuracy of gradient approximations in both the OtP and PtO approaches.
Current efforts are focused on extending the error estimates derived here to account for the time discretization. We believe this work opens several promising direction for future research, including rank-adaptive strategies based on the current optimization progress, high-order integrators for adjoint covariant ODEs, alternative discretize-then-project formulations and the extension of these techniques to complex, real-world applications.

\section*{Acknowledgments} TV is member of the INdAM-GNCS group and has been partially funded by Progetto di Ricerca GNCS-INdAM, CUP E53C25002010001. The authors acknowledge the use of large language models (ChatGPT \& Claude) limited to occasionally improving the fluency of the text and to generating part of the documentation of the repository associated with this manuscript. The authors assume responsibility for all content.

\bibliographystyle{plain}
\bibliography{references}

\appendix

\section{Adjoint of the map $H\mapsto D\cP_Y[H]Z$}\label{appendix:adjoint}
In this appendix we detail the calculations to derive the adjoint of the map $H\mapsto D\cP_Y[H]Z$ for a given $Z\in \mathbb{R}^{m\times n}$. To this end, we first consider the first contribution in \eqref{eq:derivative_projection_compact}. Let $(L,U_p,V_p)$ be a triplet identifying $H\in T_Y\cMr$. Then, for any $W\in\mathbb{R}^{m\times n}$ and due to the cyclic properties of the trace operator
\begin{equation*}
\begin{aligned}
\left\langle \left( U S^{-\top} U_p\tr + U_p S^{-1}U\tr \right)M Z (P_V^\perp)\tr, W\right\rangle_{M,\Xi}&=\trace\left(\left(\left(U S^{-\top} U_p\tr + U_p S^{-1}U\tr\right)M Z (P_V^\perp)\tr\right)\tr M W \Xi\right)\\
&=\trace\left(P_V^\perp Z\tr M \left(U_pS^{-1}U\tr + US^{-\top} U_p\tr\right) M W\Xi\right)\\
&=\trace\left(\left(U_pS^{-1}U\tr + US^{-\top} U_p\tr\right) M W\Xi P_V^\perp Z\tr M\right).
\end{aligned}
\end{equation*}
We now express $U_p$ in terms of $H$ using \eqref{eq:tangent_parametrization}, obtaining
\begin{equation*}
\begin{aligned}
\left\langle \left( U_p S^{-1}U\tr + U S^{-\top} U_p\tr \right)M Z (P_V^\perp)\tr, W\right\rangle_{M,\Xi}&=\trace\left(P_U^\perp H \Xi V S^{-1}U\tr M W\Xi P_V^\perp Z\tr M\right)\\
&+\trace\left(US^{-\top} V\tr \Xi H\tr (P_U^\perp)\tr M W\Xi P_V\tr Z\tr M\right)\\
&=\trace\left(H\tr (P_U^\perp)\tr M Z (P_V^\perp)\tr \Xi W\tr M U S^{-\top} V\tr \Xi\right)\\
&+\trace\left( H\tr (P_U^\perp)\tr M W\Xi P_V^\perp Z\tr M US^{-\top} V\tr \Xi \right)\\
&=\left\langle H, \cL_{Z,1}^\star (W)\right\rangle_{M,\Xi},
\end{aligned}
\end{equation*}
with 
\[\cL_{Z,1}^\star(W):=M^{-1} (P_U^\perp)\tr M \left[W \Xi P_V^\perp Z\tr + Z(P_V^\perp)\tr \Xi W\tr\right]M U S^{-\top}V\tr.\]
Next, we focus on the second contribution,
\begin{equation*}
\begin{aligned}
\left\langle P_U^\perp Z \Xi \left( V_p S^{-\top}V\tr + V S^{-1} V_p\tr \right), W\right\rangle_{M,\Xi}&=\trace\left(\left(P_U^\perp Z \Xi \left( V_p S^{-\top}V\tr + V S^{-1} V_p\tr \right)\right)\tr M W \Xi \right)\\
&=\trace\left(\left(V_p S^{-\top}V\tr + VS^{-1}V_p\tr\right)\Xi Z\tr (P_U^\perp)\tr M W \Xi \right).
\end{aligned}
\end{equation*}
Expressing $V_p$ in terms of $H$ using again \eqref{eq:tangent_parametrization} leads to
\begin{equation*} 
\begin{aligned}
\left\langle P_U^\perp Z \Xi \left( V_p S^{-\top}V\tr + V S^{-1} V_p\tr \right), W\right\rangle_{M,\Xi}&=\trace\left(\left(P_U^\perp Z \Xi \left( V_p S^{-\top}V\tr + V S^{-1} V_p\tr \right)\right)\tr M W \Xi \right)\\
&=\trace\left(P_V^\perp H\tr M U S^{-\top}V\tr \Xi Z\tr (P_U^\perp)\tr M W \Xi \right)\\
&+\trace\left(VS^{-1} U\tr M H (P_V^\perp)\tr \Xi Z\tr (P_U^\perp)\tr M W \Xi \right)\\
&=\trace\left(H\tr M U S^{-\top}V\tr \Xi Z\tr (P_U^\perp)\tr M W \Xi P_V^\perp \right)\\
&+\trace\left(H\tr M U S^{-\top} V\tr\Xi W\tr M P_U^\perp Z \Xi P_V^\perp  \right)\\
&=\left\langle H, \cL^\star_{Z,2}(W)\right\rangle_{M,\Xi},
\end{aligned}
\end{equation*}
where
\[\cL^\star_{Z,2}(W)=US^{-\top}V\tr\Xi \left(Z\tr (P_U^\perp)\tr M W + W\tr M P_U^\perp Z\right)\Xi P_V^\perp \Xi^{-1}.\]
We thus conclude that for every $W\in\mathbb{R}^{m\times n}$,
\[\left\langle D\cP_Y[H]Z, W\right\rangle_{M,\Xi}=\left\langle H, \cL^\star_{Z,1}(W)+\cL^\star_{Z,2}(W)\right\rangle_{M,\Xi}=:\left\langle H, D\cP^\star_{Y,Z}(W)\right\rangle_{M,\Xi}.\]
The expression of the adjoint simplifies if we restrict $W$ to be on the tangent space to $\cMr$ at $Y$. In such cases
we have that $(P_U^\perp)\tr M W \Xi P_V^\perp=0$, so that for all $W\in T_{Y}\cMr$,
\begin{equation}\label{eq:simplified_expression_adjoint}
D\cP^\star_{Y,Z}(W)= M^{-1} (P_U^\perp)\tr M Z (P_V^\perp)\tr \Xi W\tr M U S^{-\top}V\tr + US^{-\top}V\tr\Xi W\tr M P_U^\perp Z \Xi P_V^\perp \Xi^{-1}.
\end{equation}
Furthermore, the identities $M P_U^\perp = (P_U^\perp)\tr M$ and $\Xi P_V^\perp = (P_V^\perp)\tr \Xi$ lead to the final expression
\begin{equation}\label{eq:simplified_expression_adjoint2}
\begin{aligned}
D\cP^\star_{Y,Z}(W)&= P_U^\perp Z \Xi P_V^\perp W\tr M U S^{-\top}V\tr + US^{-\top}V\tr\Xi W\tr (P_U^\perp)\tr M Z (P_V^\perp)\tr\\
&=P_U^\perp Z \Xi \widetilde{V}_p S^{-\top}V\tr + US^{-\top} \widetilde{U}_p\tr M Z (P_V^\perp)\tr,
\end{aligned}
\end{equation}
where $(\widetilde{L},\widetilde{U}_p,\widetilde{V}_p)$ is the triplet that identifies $W\in T_Y\cMr$.
Note that, for any $Z\in \Rmn$ and $W\in T_{Y}\cMr$, $D\cP^\star_{Y,Z}(W)\in T_{Y}\cMr$ with the parametrization
\begin{equation}\label{eq:tangent_parametrization_adjoint}
\widehat{L}:=0,\quad \widehat{U}_p:=P_U^\perp Z \Xi \widetilde{V}_p S^{-\top},\quad \widehat{V}_p:= P_V^\perp Z\tr M \widetilde{U}_p S^{-1}.
\end{equation}
An additional simplification arises if $Z\in N_{Y}\cMr$. Letting $\widetilde{Z}\in \R^{r\times r}$ be such that $Z=P_{U}^\perp \widetilde{Z}(P_V^\perp)\tr$ it holds,
\begin{equation}\label{eq:simplified_expression_adjoint3}
\begin{aligned}
D\cP^\star_{Y,Z}(W)&=Z \Xi \widetilde{V}_p S^{-\top}V\tr + US^{-\top} \widetilde{U}_p\tr M Z\\
&=P_U^\perp \widetilde{Z} \Xi \widetilde{V}_p S^{-\top}V\tr + US^{-\top} \widetilde{U}_p\tr M \widetilde{Z} (P_V^\perp)\tr.
\end{aligned}
\end{equation}

\end{document}